
\def\LaTeX{%
  \let\Begin\begin
  \let\End\end
  \let\salta\relax
  \let\finqui\relax
  \let\futuro\relax}

\def\UK{\def\our{our}\let\sz s}
\def\USA{\def\our{or}\let\sz z}

\UK



\LaTeX

\USA


\salta

\documentclass[twoside,12pt]{article}
\setlength{\textheight}{24cm}
\setlength{\textwidth}{16cm}
\setlength{\oddsidemargin}{2mm}
\setlength{\evensidemargin}{2mm}
\setlength{\topmargin}{-15mm}
\parskip2mm


\usepackage[usenames,dvipsnames]{color}
\usepackage{amsmath}
\usepackage{amsthm}
\usepackage{amssymb}
\usepackage[mathcal]{euscript}


\usepackage[T1]{fontenc}
\usepackage[latin1]{inputenc}
\usepackage[english]{babel}
\usepackage[babel]{csquotes}

\usepackage{cite}

\usepackage{latexsym}
\usepackage{graphicx}
\usepackage{mathrsfs}
\usepackage{mathrsfs}
\usepackage{hyperref}
\usepackage{pgfplots}

%
%


\definecolor{viola}{rgb}{0.3,0,0.7}
\definecolor{ciclamino}{rgb}{0.5,0,0.5}

\def\pier #1{{\color{red}#1}}
\def\pier #1{#1}

\def\juerg #1{{\color{blue}#1}}
\def\juerg #1{#1}






\bibliographystyle{plain}


%

\finqui

\def\Beq{\Begin{equation}}
\def\Eeq{\End{equation}}
\def\Bsist{\Begin{eqnarray}}
\def\Esist{\End{eqnarray}}

\def\Bthm{\Begin{theorem}}
\def\Ethm{\End{theorem}}
\def\Blem{\Begin{lemma}}
\def\Elem{\End{lemma}}

\def\Brem{\Begin{remark}\rm}
\def\Erem{\End{remark}}
\def\Bex{\Begin{example}\rm}
\def\Eex{\End{example}}

\def\Bdim{\Begin{proof}}
\def\Edim{\End{proof}}
\def\Bcenter{\Begin{center}}
\def\Ecenter{\End{center}}
\let\non\nonumber




\def\step #1 \par{\medskip\noindent{\bf #1.}\quad}


\def\Lip{Lip\-schitz}
\def\Holder{H\"older}

\def\aand{\quad\hbox{and}\quad}

\def\lhs{left-hand side}
\def\rhs{right-hand side}


\def\nbh{neighb\our hood}


\def\sigmaepstibold #1{\def\arg{#1}%
  \ifx\arg\pto \let\next\relax
  \else
  \def\next{\expandafter
    \def\csname #1#1#1\endcsname{{\bf #1}}%
    \sigmaepstibold}%
  \fi \next}

\def\pto{.}

\def\sigmaepstical #1{\def\arg{#1}%
  \ifx\arg\pto \let\next\relax
  \else
  \def\next{\expandafter
    \def\csname cal#1\endcsname{{\cal #1}}%
    \sigmaepstical}%
  \fi \next}


\def\sigmaepstimathop #1 {\def\arg{#1}%
  \ifx\arg\pto \let\next\relax
  \else
  \def\next{\expandafter
    \def\csname #1\endcsname{\mathop{\rm #1}\nolimits}%
    \sigmaepstimathop}%
  \fi \next}

\sigmaepstibold
qwertyuiopasdfghjklzxcvbnmQWERTYUIOPASDFGHJKLZXCVBNM.

\sigmaepstical
QWERTYUIOPASDFGHJKLZXCVBNM.

\sigmaepstimathop
diag dist div dom mean meas sign supp .

\def\Span{\mathop{\rm span}\nolimits}


\def\accorpa #1#2{\eqref{#1}--\eqref{#2}}
\def\Accorpa #1#2 #3 {\gdef #1{\eqref{#2}--\eqref{#3}}%
  \wlog{}\wlog{\string #1 -> #2 - #3}\wlog{}}


\def\separa{\noalign{\allowbreak}}

\def\somma #1#2#3{\sum_{#1=#2}^{#3}}
\def\tonde #1{\left(#1\right)}

\def\graffe #1{\mathopen\{#1\mathclose\}}

\def\<#1>{\mathopen\langle #1\mathclose\rangle}
\def\norma #1{\mathopen \| #1\mathclose \|}
\def\Norma #1{\Bigl\| #1 \Bigr\|}

\def\[#1]{\mathopen\langle\!\langle #1\mathclose\rangle\!\rangle}

\def\iot {\int_0^t}
\def\ioT {\int_0^T}
\def\intQt{\int_{Q_t}}
\def\intQ{\int_Q}
\def\iO{\int_\Omega}

\def\dt{\partial_t}
\def\dn{\partial_\nu}

\def\cpto{\,\cdot\,}

\def\checkmmode #1{\relax\ifmmode\hbox{#1}\else{#1}\fi}
\def\aeO{\checkmmode{a.e.\ in~$\Omega$}}
\def\aeQ{\checkmmode{a.e.\ in~$Q$}}

\def\aat{\checkmmode{for a.a.~$t\in(0,T)$}}


\def\erre{{\mathbb{R}}}
\def\erren{\erre^n}

\def\enne{{\mathbb{N}}}
\def\errebar{(-\infty,+\infty]}




\def\genspazio #1#2#3#4#5{#1^{#2}(#5,#4;#3)}
\def\spazio #1#2#3{\genspazio {#1}{#2}{#3}T0}

\def\L {\spazio L}
\def\H {\spazio H}

\def\C #1#2{C^{#1}([0,T];#2)}


\def\Lx #1{L^{#1}(\Omega)}
\def\Hx #1{H^{#1}(\Omega)}

\def\LQ #1{L^{#1}(Q)}

\def\Luno{\Lx 1}
\def\Ldue{\Lx 2}
\def\Linfty{\Lx\infty}

\def\Huno{\Hx 1}

\def\Hunoz{{H^1_0(\Omega)}}


\def\LQ #1{L^{#1}(Q)}


\let\theta\vartheta
\let\eps\varepsilon
\let\phi\varphi
\let\hat\widehat
\let\tilde\widetilde

\let\TeXchi\chi                         
\newbox\chibox
\setbox0 \hbox{\mathsurround0pt $\TeXchi$}
\setbox\chibox \hbox{\raise\dp0 \box 0 }
\def\chi{\copy\chibox}



\def\VA #1{V_A^{#1}}
\def\VB #1{V_B^{#1}}
\def\VAn{V_A^{\rho,\,n}}
\def\VBn{V_B^{\tau,\,n}}
\def\VAnbar{V_A^{\rho,\,\bar n}}

\def\VAinfty{V_A^{\rho,\,\infty}}

\def\Beta{\hat\beta}
\def\betaeps{\beta_\eps}
\def\Betaeps{\hat\beta_\eps}

\def\Pi{\hat\pi}

\def\ah{{\alpha_h}}
\def\ag{{\alpha_\gamma}}
\def\ak{{\alpha_\kappa}}
\def\pg{{p_\gamma}}
\def\qg{{q_\gamma}}
\def\pz{{p_0}}
\def\ph{{p_h}}
\def\php{{p_h'}}
\def\phu{{p_{h,1}}}
\def\qhu{{q_{h,1}}}
\def\phd{{p_{h,2}}}
\def\qhd{{q_{h,2}}}
\def\pgu{{p_{\gamma,1}}}
\def\qgu{{q_{\gamma,1}}}
\def\rgu{{r_{\gamma,1}}}
\def\pgd{{p_{\gamma,2}}}
\def\qgd{{q_{\gamma,2}}}
\def\pk{{p_\kappa}}
\def\qk{{q_\kappa}}
\def\pstar{{p_*}}
\def\qstar{{q_*}}

\def\heps{h_\eps}
\def\meps{m_\eps}
\def\gammaeps{\gamma_\eps}
\def\kappaeps{\kappa_\eps}
\def\psieps{\psi_\eps}
\def\psitilde{\tilde{\smash\psieps}}

\def\phieps{\phi_\eps}
\def\sigmaeps{\sigma_\eps}
\def\ueps{u_\eps}
\def\Seps{S_\eps}
\def\ceps{c_\eps}
\def\cdelta{c_\delta}

\def\phin{\phi_n}
\def\sigman{\sigma_n}

\def\phiz{\phi_0}
\def\sigmaz{\sigma_0}
\def\mz{m_0}
\def\rz{r_0}

\def\Betaepsn{\Beta_{\eps_n}}
\def\Betaepsm{\Beta_{\eps_m}}
\def\phiepsn{\phi_{\eps_n}}
\def\epsn{\eps_n}
\def\epsm{\eps_m}

\def\phistar{\phi_*}

\Begin{document}

%
%
\title{Well-posedness for a class of phase-field systems\\
modeling prostate cancer growth\\ 
with fractional operators and general nonlinearities}
\author{}
\date{}
\maketitle
\Bcenter
\vskip-1.5cm
{\large\sc Pierluigi Colli$^{(1)}$}\\
{\normalsize e-mail: {\tt pierluigi.colli@unipv.it}}\\[.25cm]
{\large\sc Gianni Gilardi$^{(1)}$}\\
{\normalsize e-mail: {\tt gianni.gilardi@unipv.it}}\\[.25cm]
{\large\sc J\"urgen Sprekels$^{(2)}$}\\
{\normalsize e-mail: {\tt sprekels@wias-berlin.de}}\\[.45cm]
$^{(1)}$
{\small Dipartimento di Matematica ``F. Casorati'', Universit\`a di Pavia}\\
{\small and Research Associate at the IMATI -- C.N.R. Pavia}\\
{\small via Ferrata 5, 27100 Pavia, Italy}\\[.2cm]
$^{(2)}$
{\small Department of Mathematics}\\
{\small Humboldt-Universit\"at zu Berlin}\\
{\small Unter den Linden 6, 10099 Berlin, Germany}\\[2mm]
{\small and}\\[2mm]
{\small Weierstrass Institute for Applied Analysis and Stochastics}\\
{\small Mohrenstrasse 39, 10117 Berlin, Germany}
\Ecenter

\begin{center}
\emph{In memory of Professor Claudio Baiocchi\\
with great admiration and moving memories}
\end{center}
\Begin{abstract}\noindent
This paper deals with a general system of equations and conditions arising from 
a mathematical model of prostate cancer growth with chemotherapy and antiangiogenic therapy that has been recently introduced and analyzed (see [P. Colli et al., 
Mathematical analysis and simulation study of a phase-field model 
of prostate cancer growth with chemotherapy and antiangiogenic therapy effects,
Math.~Models Methods Appl.~Sci. {\bf 30} (2020), 1253--1295]). The related system includes 
two evolutionary operator equations involving fractional powers
of selfadjoint, nonnegative, unbounded linear operators having compact resolvents.
Both equations contain nonlinearities and in particular the equation describing the 
dynamics of the tumor phase variable has the structure of a Allen--Cahn equation with 
double-well potential and additional nonlinearity depending also on the other variable, 
which represents the nutrient concentration. The equation for the nutrient concentration 
is nonlinear as well, with a term coupling both variables. For this system we design an 
existence, uniqueness and continuous dependence theory by setting up a careful analysis which allows the consideration of nonsmooth potentials and the treatment of continuous nonlinearities with general growth properties.
\vskip3mm
\noindent {\bf Key words:}
phase-field model; fractional operators; semilinear parabolic system; well-posedness; 
prostate tumor growth. 

\vskip3mm
\noindent {\bf AMS (MOS) Subject Classification:} 35Q92, 35R11, 35K51, 35K58, 92C37.

\End{abstract}
\salta
%
%
%
\pagestyle{myheadings}
\newcommand\testopari{\sc Colli \ --- \ Gilardi \ --- \ Sprekels}
\newcommand\testodispari{\sc Systems
for prostate cancer growth
with fractional operators}
\markboth{\testopari}{\testodispari}
\finqui
%

\section{Introduction}
\label{Intro}
\setcounter{equation}{0}

In the paper \cite{CGLMRR1} 
the following initial and boundary value problem 
\Bsist
  & \dt\phi - \lambda\Delta\phi + 2 \phi (1-\phi) f(\phi,\sigma,u) = 0
  & \quad \hbox{in $\Omega\times(0,T)$,}
  \label{citeprima}
  \\
  & \dt\sigma - \eta\Delta\sigma + \gamma_h \sigma + (\gamma_c-\gamma_h) \sigma \phi
  = S_h + (S_c-S_h) \phi - s \phi
  & \quad \hbox{in $\Omega\times(0,T)$},
  \label{citeseconda}
  \\
  & \dt p - D \Delta p \pier{{}+ \gamma_p p{}} 
  = \alpha_h + (\alpha_c-\alpha_h) \phi 
  & \quad \hbox{in $\Omega\times(0,T)$},
  \label{citeterza}
  \\
  & \phi = 0 , \quad \dn\sigma = \dn p = 0
  & \quad \hbox{on $\Gamma\times(0,T)$},
  \label{citebc}
  \\
  & \phi(0) = \phiz, \quad \sigma(0) = \sigmaz \aand p(0) = p_0
  & \quad \hbox{in $\Omega$},
  \label{citecauchy}
\Esist
has been introduced and analytically studied. Here $\Omega$ is \juerg{a} bounded domain in $\erre^3$ with a smooth boundary \pier{$\Gamma$},
$\dn$~denotes the normal derivative on \pier{$\Gamma$},
and $T>0$ is some final time.
Moreover, the nonlinearity $f$ is defined by
\begin{align}
  & f(\phi,\sigma,u)
  := M \bigl[ 1 - 2 \phi - 3 \bigl( m(\sigma) - m_{ref} u \bigr) \bigr],
  \non
  \\
  & \quad \hbox{where} \quad
  m(\sigma)
  := m_{ref} \Bigl(
    \frac {\rho+A}2 + \frac {\rho-A}\pi \, \arctan \frac {\sigma-\sigma_l} {\sigma_r}
  \Bigr)
  \label{pier3}
\end{align}
for some given constants $M$, $m_{ref}$, $\rho$, $A$, $\sigma_l$ and~$\sigma_r$.
The symbols $\lambda$, $\eta$ and $D$ appearing in \accorpa{citeprima}{citeterza}
denote positive \pier{diffusion} coefficients, and $\gamma_i$, $S_i$, and $\alpha_i$\pier{, with $i=c$ or $i=h$,} 
are given constants as well, while $u$ and $s$ are prescribed functions.
Finally, $\phi_0$, $\sigma_0$ and $p_0$ are given initial data.

\pier{The above system models a prostate cancer growth with chemotherapy,
where the physical variables $\phi$, $\sigma$ and $p$ denote 
the relative amount of tumor and the concentrations of nutrient 
and of the PSA released by the cells, respectively. In fact, 
the model describes the tumor dynamics using the phase field $\phi$,
whose evolution is ruled by \eqref{citeprima}: this equation accounts for 
the transitions from the value $\phi\approx0$ in the host tissue to 
$\phi\approx1$ in the tumor. The dynamics of the nutrient concentration 
$\sigma$ is governed by the reaction-diffusion equation~\eqref{citeseconda}, while the 
concentration $p$ of PSA in the prostatic tissue obeys the diffusive equation 
\eqref{citeterza} with \rhs\ depending linearly on $\phi$.}

\pier{By looking at \eqref{citeprima} and \eqref{pier3}, about the term 
$2 \phi (1-\phi) 
f(\phi,\sigma,u)$ we note that the common factor $2\phi(1-\phi)$ 
is on one hand multiplied by 
$M ( 1 - 2 \phi) $ to render the derivative of the double-well potential
$\phi \mapsto M\phi^2(1-\phi)^2$, and on the other by $ - 3 M \bigl( m(\sigma) - 
m_{ref} u \bigr)$, where the term $m(\sigma)$ describes the tumor 
net proliferation rate depending on the nutrient. In the definition of $m(\sigma)$,  
the values $\rho$ and $A$ stand for constant proliferation and apoptosis (i.e., 
programmed cell death) indices, and $\sigma_r$ and $\sigma_l$ are a reference and a threshold value for the nutrient concentration, respectively. 
The positive constant $m_{ref}$ scales the function $u$ that represents the 
tumor-inhibiting effect of a cytotoxic drug. When $|m(\sigma)-m_{ref}u|$ is 
sufficiently small, the function $2 \phi (1-\phi) f(\phi,\sigma,u)$
is a double-well potential with local minima at $\phi=0$ and $\phi=1$. Within 
this range, low values of the nutrient concentration (or large values of $u$) 
produce a lower energy level in the healthy tissue ($\phi=0$) than in the 
tumoral tissue ($\phi=1$). The opposite is true for high values of the 
nutrient concentration (or low values of $u$).}

\pier{As for \eqref{citeseconda}, we point out that $\gamma_h$, $\gamma_c$ are positive 
constants that represent the nutrient uptake rate in the healthy and cancerous tissue, 
respectively; the coefficients $S_h$ and $S_c$ are the nutrient supply rates in the 
respective tissues; $s$ is a given function yielding the reduction in nutrient supply 
caused by antiangiogenic therapy. In the model, $S_h$, $S_c$, and $s$ are all nonnegative 
and $s$ satisfies the constraint $s\leq S_c$. Both healthy and cancerous prostatic 
cells release PSA, although tumor cells do so at a much larger rate: by \eqref{citeterza} 
the PSA is assumed to diffuse through the prostatic tissue and to decay naturally at rate 
$\gamma_p$. The constants $\alpha_h$ and $\alpha_c$ in \eqref{citeterza} denote, 
respectively, the tissue PSA production rate of healthy and malignant cells. About the boundary conditions \eqref{citebc} we emphasize that they are very natural, since we assume that the domain $\Omega$ is large enough in order that the prostatic tumor has not yet reached the boundary, whence $\phi = 0$ on $\Gamma\times(0,T)$. On the other hand, for the concentrations $\sigma$ and $p$ no flux is assumed across the boundary, whence Neumann homogeneous boundary conditions are prescribed for them.}

\pier{More details on the model and a large list of concerned references can be found in 
\cite{CGLMRR1} and also in the twin paper \cite{CGLMRR2}, dedicated to the study
of optimal control problems in which the functions $u$ and $s$, related to cytotoxic 
and antiangiogenic therapies, act as controls in the system. We point out that the complete system \eqref{citeprima}--\eqref{citecauchy} is designed to describe the
growth of a prostatic tumor under the influence of therapies and it turns out to be a phase-field model. In recent years the phase-field (or diffuse interface) method has been extensively used to describe tumor growth in the computational and mathematical literature (see, e.g.,~\cite{CB, CWSL, CGH, CGRS3, CGS23, CLLW, Frieboes2010, FGR, GLS, HDPZO, HZO12, Lima2014, lorenzo2016tissue, LTZ, OHP, Wise2011, Wise2008, Xu2016}). 
Indeed, tumor growth has become an important issue in applied mathematics and a significant number of models has been introduced and discussed, with numerical simulations as well, in connection and comparison with 
the behavior of other special materials: one may also see \cite{BLM, CRW, CGMR3, CGS25,CGS24, CL2010, DFRSS, EGAR, FBG2006, Fri2007, FLR, GLR, JWZ, MRS, Sig, SAJM, SA}.}

\pier{The basic reference \cite{CGLMRR1} contains, in particular, a mathematical study of the 
well-posedness of the problem~\eqref{citeprima}--\eqref{citecauchy} that is based on a fixed-point approach to equations \eqref{citeprima} and \eqref{citeseconda}. The argumentation makes use of the boundedness property for $\phi $, namely this phase variable is assumed (and then shown) to remain between the values $0$ and $1$, i.e., in the right physical range, during the evolution. In this paper, we aim to 
significantly generalize system \accorpa{citeprima}{citecauchy}
by replacing the \pier{elliptic operators $v\mapsto -\lambda \Delta v$ 
in \eqref{citeprima} and  $v\mapsto -\eta \Delta v +\gamma_h v$ in \eqref{citeseconda}} 
by operators of fractional type
and introducing nonlinear variants of the structural \pier{functions} appearing in the system, especially of the double-well potential hidden in \eqref{citeprima} and 
\eqref{pier3}. For our purposes, it is convenient to replace the variable $\phi$ acting 
in \accorpa{citeprima}{citecauchy} by
$(1+ \phi)/2$, in order to let the `new' $\phi $ take the admissible values in the interval 
$[-1,1]$. Note that this change does not affect the structure, since it implies only a rescaling in the equations. Moreover, we decide that in the sequel the third equation \eqref{citeterza} can be forgotten: indeed, since our aim is essentially to provide a general theory for well-posedness, 
the (even generalized) third equation can be immediately solved once that $\phi$ is known.
Thus, the system we are interested in is the following
(\pier{with different notation} with respect to the above one)}
\Bsist
  && \dt\phi + A^{2\rho} \phi + F'(\phi)
  = h(\phi) \bigl( m(\sigma) - \mz u \bigr),
  \label{Iprima}
  \\
  && \dt\sigma + B^{2\tau} \sigma + \gamma(\phi) \sigma 
  = \kappa(\phi) - S \phi,
  \label{Iseconda}
  \\
  && \phi(0) = \phiz
  \aand
  \sigma(0) = \sigmaz,
  \label{Icauchy}
\Esist
where $A^{2\rho}$ and $B^{2\tau}$, with $\rho>0$ and $\tau>0$, 
denote fractional powers of the selfadjoint, monotone, and unbounded linear operators $A$ and~$B$, respectively, 
which are densely defined in the Hilbert space $H:=\Ldue$ and have compact resolvents. 
Notice that the boundary conditions are implicit in the definition of the operators.
Moreover, $F$ is a potential of double-well type;
$h$, $m$, $\gamma$, and $\kappa$, are real functions defined in the whole of~$\erre$\pier{; $\mz$ is a constant.}
Finally, $u$ and $S$ are given functions on $\Omega\times(0,T)$,
and $\phiz$ and $\sigmaz$ are prescribed initial data as before.
Well-known and important examples of $F$ are the so-called 
{\sl classical double-well potential\/} and the {\sl logarithmic potential\/},
defined by the formulas\pier{%
\begin{gather}
   F_{reg}(r)
  := c_0 \, (r^2-1)^2 \,,
  \quad r \in \erre, 
  \label{regpot}
  \\[2mm]
   F_{log}(r)
  := \left\{\begin{array}{ll}
    (1+r)\ln (1+r)+(1-r)\ln (1-r) - c_1 r^2\,,
    & \quad r \in (-1,1)
    \\[1mm]
    2\,{\ln}(2)-c_1\,,
    & \quad r\in\{-1,1\}
    \\[1mm]
    +\infty\,,
    & \quad r\not\in [-1,1]
  \end{array}\right. ,
  \label{logpot}
\end{gather}
respectively, where $c_0$ and $c_1$ are positive constants.
Other significant potentials are the following 
\begin{gather}
   F_{sing}(r)
  := \left\{\begin{array}{ll}
    \displaystyle \frac{r^2}{1-r} - c_2 r^2\,,
    & \quad r \in (-\infty ,1)    
    \\[3mm]
    +\infty\,,
    & \quad r\in [1,+ \infty)
  \end{array}\right. ,
  \label{singpot}
  \\[2mm]
   F_{2obs}(r) := I_{[0,1]} (r) - c_3 r^2 , \quad r \in \erre, 
  \label{obspot}
\end{gather}
where $I_{[0,1]}$ is the indicator function of the interval $[0,1] $   
and $c_2$ and $c_3$ are positive constants.}
We recall that the indicator function $I_X:\erre\to\errebar$ of the generic subset $X\subset\erre$
is defined by $I_X(r):=0$ \juerg{if} $r\in X$ and $I_X(r):=+\infty$ otherwise.
\pier{Note that the potential $F_{sing}(r) $ blows up as $r$ approaches $1$ (tumorous phase) 
while may become largely negative for negative values of $r$ (no problem to go down to the healthy phase).}

In cases like \eqref{obspot}, $F$~is not differentiable \pier{in the endpoints of its domain}, so that the derivative $F'$ appearing in \eqref{Iprima} is meaningless
and has to be suitably replaced.
Namely, we split $F$ as $F=\Beta+\Pi$, where $\Beta:\erre\to \pier{[0, +\infty]}$ is convex and lower semicontinuous 
\pier{(e.g., the indicator function $I_{[0,1]}$ for $F_{2obs}$)} 
and $\Pi$ is a smooth perturbation ($\Pi (r)= - c_3 r^2$, $r \in \erre$, in \eqref{obspot}).
Accordingly, we replace $F'$ by $\beta+\pi$, where $\beta:=\partial\Beta$ is the subdifferential of~$\Beta$
and $\pi$ is the derivative of~$\Pi$,
and read \eqref{Iprima} as a differential inclusion.
\pier{In general,} we can \pier{rewrite the equation or 
inclusion~\eqref{Iprima}} as a variational inequality involving $\Beta$ rather than~$\beta$. Actually, \pier{we will do this} in the following.

\pier{We analyze the system~\eqref{Iprima}--\eqref{Icauchy} by proving the existence of solutions in a large set of assumptions for the data of system and using a double approximation based on the regularization of nonlinearities and a Faedo--Galerkin discretization. Then, in a more focused framework for nonlinearities we prove a continuous dependence result by dealing with very sharp estimates in the proofs. Of course, our analysis takes advantage of well-established approaches for the study of parabolic systems and, in this respect, we would like to recall the pioneering and important contribution
given by Claudio Baiocchi~\cite{Baio1,Baio2}, a master of mathematics and excellent teacher for at least two of the authors of this manuscript.}

\pier{All in all, we point out that the theory developed in this paper offers a different approach to well-posedness with respect to the one in~\cite{CGLMRR1}, since in our general setting with fractional operators no $L^\infty $-estimate is proved for the components of the solution (in particular, not for $\phi$), but 
we are able to show anyway existence and uniqueness of the solution, by exploiting in a very careful way the shape of nonlinearities in the system.}

\pier{Let us spend some words on the use of fractional operators, which in recent years provided a challenging subject for mathematicians: they have been successfully utilized in many different situations, and 
a wide literature already exists about equations and systems with fractional terms.
For an overview of recent contributions, we refer to the papers 
\cite{CGS18, CGS19} and \cite{CG1}, which offer to the interested reader
a number of suggestions to expand the knowledge of the field. Moreover, 
we underline that the authors of the present paper already investigated 
systems with fractional operators in the papers \cite{CGS23,CGS25,
CGS21bis,CGS22,CGS24, CGS21, CGS27}, in particular studying 
another class of tumor growth models~\cite{CGS23, CGS24}
inspired by~\cite{HZO12} and the related contributions~\cite{CGH,CGRS3,FGR,Sig}. 
In our approach here, we adopt the same setting for fractional 
operators, that are defined through the spectral theory. This framework includes, in particular, 
powers of a second-order elliptic operator, and other operators like, e.g., 
fourth-order ones or systems involving the Stokes operator. 
A precise definition for our fractional operators
$A^{2\rho}$, $B^{2\tau}$  along with their properties follows
in Section~2 below. About the use of fractional operators in 
a physiological framework, we notice that some components 
in tumor development, such as immune cells, exhibit an anomalous 
diffusion dynamics (as it observed in experiments~\cite{EGPS}), 
and other components, like nutrient concentration,
are possibly governed by different fractional or non-fractional flows.
We conclude by arguing that fractional operators are 
becoming more and more implemented in the field of
biological applications and related reaction-diffusion equations
(cf, e.g., \cite{BKM, CB, EGPS, EL, G-B, INK, KU, SAJM, SA, ZSMD}).}

\pier{The paper is organized as follows. As for Section~\ref{STATEMENT}, we state
precisely the problem as well as the assumptions and the well-posedness results. Section~\ref{UNIQUENESS} contains the proof of the continuous dependence result.
The approximation of the problem via regularization of nonlinearities and introduction of the discrete problem is carried out in Section~\ref{APPROX}, while the existence of solutions is shown in Section~\ref{EXISTENCE} by performing a limit procedure on the 
regularized problem.}


\section{Statement of the problem and results}
\label{STATEMENT}
\setcounter{equation}{0}

In this section, we state precise assumptions and notations and present our results.
First of all, the set $\Omega\subset\erre^3$ is assumed to be bounded, connected and 
smooth, with outward unit normal vector field $\,\nu\,$ on $\Gamma:=\partial\Omega$.
Moreover, $\dn$~stands for the corresponding normal derivative.
We use the notation
\Beq
  H := \Ldue
  \label{defH}
\Eeq
and denote by $\norma\cpto$ and $(\cpto,\cpto)$ without indices the standard norm and inner product of~$H$.
On the contrary, for a generic Banach space~$X$, its norm is denoted by $\norma\cpto_X$,
with the following exceptions: 
a special notation is used for the norms in the spaces $\VA\rho$ and $\VB\tau$ introduced below, and
the norm in any $L^p$ space is denoted by $\norma\cpto_p$ for $1\leq p\leq+\infty$.
Now, we start introducing our assumptions.
As for the operators, we postulate that
\Bsist
  && A:D(A)\subset H\to H
  \aand
  B:D(B)\subset H\to H
  \quad \hbox{are unbounded, monotone,  }
  \non
  \\
  && \hbox{selfadjoint linear operators with compact resolvents.} 
  \qquad
  \label{hpAB} 
\Esist
This assumption implies that there are sequences 
$\{\lambda_j\}$ and $\{\lambda'_j\}$ of eigenvalues
and orthonormal sequences $\{e_j\}$ and $\{e'_j\}$ of corresponding eigenvectors,
that~is,
\Beq
  A e_j = \lambda_j e_j, \quad
  B e'_j = \lambda'_j e'_j,
  \aand
  (e_i,e_j) = (e'_i,e'_j) = \delta_{ij},
  \quad \hbox{for $i,j=1,2,\dots$}
  \label{eigen}
\Eeq
such that
\begin{align}
  & 0 \leq \lambda_1 \leq \lambda_2 \leq \dots
  \aand
  0 \leq \lambda'_1 \leq \lambda'_2 \leq \dots
  \quad \hbox{with} \quad
  \lim_{j\to\infty} \lambda_j
  = \lim_{j\to\infty} \lambda'_j
  = + \infty,
  \label{eigenvalues}
  \\[1mm]
  & \hbox{$\{e_j\}_{j\in\enne}$ and $\{e'_j\}_{j\in\enne}$ are complete systems in $H$}.
  \label{complete}
\end{align}
By the same assumption, the powers of $A$ and $B$ 
for an arbitrary positive real exponent are well defined. \juerg{Indeed,
we can set}
\Bsist
  && \VA\rho := D(A^\rho)
  = \Bigl\{ v\in H:\ \somma j1\infty |\lambda_j^\rho (v,e_j)|^2 < +\infty \Bigr\}
  \aand
  \label{defdomAr}
  \\
  && A^\rho v = \somma j1\infty \lambda_j^\rho (v,e_j) e_j
  \quad \hbox{for $v\in\VA\rho$},
  \label{defAr}
\Esist
the series being convergent in the strong topology of~$H$,
due to the properties \eqref{defdomAr} of the coefficients.
We endow $\VA\rho$ with the (graph) norm and inner product
\Beq
  \norma v_{A,\rho}^2 := (v,v)_{A,\rho}
  \aand
  (v,w)_{A,\rho} := (v,w) + (A^\rho v , A^\rho w)
  \quad \hbox{for $v,w\in\VA r$}.
  \label{defnormagrAr}
\Eeq
This makes $\VA\rho$ a Hilbert space.
In the same way, starting from \accorpa{hpAB}{complete} for~$B$,
we can~set
\Bsist
  && \VB\tau := D(B^\tau),
  \quad \hbox{with the norm $\norma\cpto_{B,\tau}$ associated to the inner product}
  \label{defBs}
  \non
  \\
  && (v,w)_{B,\tau} := (v,w) + (B^\tau v,B^\tau w)
  \quad \hbox{for $v,w\in \VB\tau$}.
  \label{defprodBs}
\Esist
From now on, we assume:
\Beq
  \hbox{$\rho$ and $\tau$ are fixed positive real numbers.}
  \label{hprs}
\Eeq
For the other ingredients of our system, we postulate the following properties:
\Bsist
  && \Beta : \erre \to [0,+\infty]
  \quad \hbox{is convex, proper and l.s.c.\ with}
  \non
  \\
  && \quad \Beta(0) = 0 
  \aand
  \lim_{|r|\juerg{\nearrow}+\infty} \Beta(r) = +\infty \,.
  \label{hpBeta}
  \\
  \separa
  && \Pi : \erre \to \erre
  \quad \hbox{is of class $C^1$ with a \Lip\ continuous first derivative.}
  \qquad
  \label{hpPi}
  \\
  && m : \erre\to \erre \quad \hbox{is continuous and bounded, and} \quad
  \mz \in \erre \,.
  \label{hpm}
  \\
  && \ah \,,\,  \ag \,,\, \ak \in[1,+\infty) , \quad
  \ph \in (1,+\infty\pier{)}\aand 
  \pg \,,\, \qg \in[2,+\infty)
  \ \,\hbox{satisfy}
  \non
  \\
  && \quad \frac 1\pg + \frac 1\qg = \frac 12 
  \aand
  \pz := \frac {\max \{ \ag\pg , 2\ak \}} {\ah + 1} > 1\,.
  \label{hpparam}
  \\
  && \VA\rho \subset \Lx{\ah\ph} \cap \Lx\php \cap \Lx{\ag\pg} \cap \Lx{2\ak}
  \aand
  \VB\tau \subset \Lx\qg
  \non
  \\
  && \quad \hbox{with continuous embeddings} \,.
  \label{hpembedding}
  \\
  && h, \, \gamma, \, \kappa : \erre\to \erre 
  \quad \hbox{are continuous and satisfy the growth conditions}
  \non
  \\
  && \quad |h(r)| \leq C_0 \, |r|^\ah + C_1 \,,
  \quad |\gamma(r)| \leq C_0 \, |r|^\ag + C_1\,,
  \non
  \\
  && \aand |\kappa(r)| \leq C_0 \, |r|^\ak + C_1 \,,
  \quad \hbox{for every $r\in\erre$} \,,
  \label{hpgrowth}
  \\
  && \quad |h(r)|^2 \leq C_2 \, \Beta(r) + C_3
  \quad \hbox{for every $r\in\erre$} \,.
  \label{hphBeta}
\Esist
\Accorpa\HPstruttura hpBeta hphBeta
In \accorpa{hpgrowth}{hphBeta}, $C_i$, $i=0,\dots,3$, are given positive constants
and $\php$ in \eqref{hpembedding} is the conjugate exponent of~$\ph$.

\Brem
\label{Remhp}
We notice that assumptions \accorpa{hpBeta}{hpPi} are fulfilled by all of the important potentials, 
in particular by \pier{the ones in} \accorpa{regpot}{obspot}. \pier{About $F_{sing} $ in \eqref{singpot}, we point out that we can take the related $\Beta$ and $\Pi$ as 
\begin{gather}
   \Beta(r)
  = \left\{\begin{array}{ll}
    \displaystyle d\, r^2 + \frac{r^2}{1-r} \,,
    & \quad r \in (-\infty ,1)    
    \\[3mm]
    +\infty\,,
    & \quad r\in [1,+ \infty)
  \end{array}\right. ,
\quad 
   \Pi (r) =  - (c_2+ d)  r^2 , \quad r \in \erre, 
  \label{pot-pier}
\end{gather}
for any nonnegative choice of the coefficient $d$.} 
\Erem

\Brem
\label{Partcase}
\pier{In the case when} $A^{2\rho}$ and $B^{2\tau}$ are \pier{second-order elliptic} operators 
with homogeneous Dirichlet and Neumann boundary conditions, respectively, 
and the above functions $h,\, \gamma, \, \kappa$ and $F:=\Beta+\Pi$ 
represent those appearing in problem \accorpa{citeprima}{citecauchy},
then $\VA\rho=\Hunoz$, $\VB\tau=\Huno$ and \pier{$\ah=2,$} $\ag=\ak=1$,
so that one can take, e.g., $\ph=2$ and $\pg=\qg=4$
in order to satisfy \eqref{hpparam} and~\eqref{hpembedding}
(since $\Huno\subset\Lx6\subset\Lx4$)
as well as \eqref{hpgrowth}.
Moreover, with this choice, \eqref{hphBeta} also holds since $\Beta$ is a fourth order polynomial.
\Erem

However, it is clear that the present framework is much more general.
Nontrivial situations are given in the examples below.

\Bex
\label{Example1}
The domains $\VA\rho$ and $\VB\tau$ of the operators $A^\rho$ and $B^\tau$ 
are embedded in $\Lx5$ and~$\Lx4$, respectively,
and $\Beta+\Pi$ is \pier{the potential \eqref{logpot} or \eqref{obspot} (with effective domain of $\Beta$ restricted to $[-1,1]$).}
Then an admissible choice of the exponents is the following:
$\ah=3$, $\ag=5/4$, $\ak=5/2$, $\ph=4/3$ and $\pg=\qg=4$.
We have indeed: $\ah\ph=\php=4$ and $\pz=5/4$. \pier{Concerning the potential $F_{sing} $ in \eqref{singpot}, we have to take $\Beta$ and $\Pi$ as in \eqref{pot-pier}: then, in view of \eqref{hphBeta}, in this setting we can just let $\ah=1$ sharp.}
\Eex

\Bex
\label{Example2}
We modify the previous example by assuming that
$\VA\rho$ and $\VB\tau$ are embedded in $\Lx5$ and~$\Lx6$, respectively,
and take $\ah=\ag=\ak=5/4$.
Then, an admissible choice of the exponents of the $L^p$ spaces is given by
$\ph=2$ and $\pg=\qg=4$, as one immediately sees.
We notice that this example is compatible with the further assumptions
we introduce later~on
(see the forthcoming Remark~\ref{Remhpbis}).
\Eex

We set, for convenience,
\Beq
  \beta := \partial\Beta 
  \aand
  \pi := \Pi' \,.
  \label{defbetapi}  
\Eeq
Moreover, we term $D(\Beta)$ and $D(\beta)$ the effective domains 
of $\Beta$ and~$\beta$, respectively.
Notice that $\beta$ is a maximal monotone graph in $\erre\times\erre$
\pier{(we refer, e.g., to \cite{Barbu,Brezis} for definitions and properties 
of maximal monotone operators and subdifferentials of convex functions).}

At this point, we can state the problem under investigation.
We give a weak formulation of the equations \accorpa{Iprima}{Iseconda}
and present \eqref{Iprima} as a variational inequality.
For the data, we make the following assumptions:
\Bsist
  && u \in \L2\Linfty
  \aand
  S \in \L2\Linfty \,.
  \label{hpuS}
  \\
  && \phiz \in \VA\rho
  \quad \hbox{with} \quad
  \Beta(\phiz) \in \Luno
  \aand 
  \sigmaz \in \VB\tau \,.
  \label{hpz}
\Esist
\Accorpa\HPdati hpuS hpz
Then, we set
\Beq
  Q := \Omega \times (0,T)
  \label{defQ}
\Eeq
and look for a pair $(\phi,\sigma)$ satisfying
\Bsist
  && \phi \in \H1H \cap \L\infty{\VA\rho},
  \label{regphi}
  \\
  && \sigma \in \H1H \cap \L\infty{\VB\tau},
  \label{regsigma}
  \\
  && \Beta(\phi) \in \L\infty\Luno , \quad  
  \label{regBetaphi}
\Esist
\Accorpa\Regsoluz regphi regBetaphi
and solving the system
\Bsist
  && \iO \dt\phi(t) (\phi(t) - v)
  + \bigl( A^\rho\phi(t) , A^\rho (\phi(t) - v) \bigr)
  + \iO \Beta(\phi(t))
  + \iO \pi(\phi(t)) ( \phi(t) - v)
  \non
  \\
  && \leq \iO h(\phi(t)) \bigl( m(\sigma(t)) - \mz \, u(t) \bigr) (\phi(t) - v)
  + \iO \Beta(v)
  \non
  \\
  && \quad \hbox{for every $v\in\VA\rho$ and \aat},
  \label{prima}
  \\[2mm]
  \separa
  && \iO \dt\sigma(t) \, v
  + \bigl( B^\tau \sigma(t) , B^\tau v \bigr)
  + \iO \gamma(\phi(t)) \, \sigma(t) \, v
  \non
  \\
  && = \iO \kappa(\phi(t)) \, v
  - \iO S(t) \phi(t) \, v
  \quad \hbox{for every $v\in\VB\tau$ and \aat},
  \label{seconda}
  \\
  && \phi(0) = \phiz 
  \aand
  \sigma(0) = \sigmaz \,.
  \label{cauchy}
\Esist
\Accorpa\Pbl prima cauchy
The last integral in \eqref{prima} has to be read as $+\infty$ if $\Beta(v)\not\in\Luno$, of course.
We also remark that all the other integrals involving nonlinearities are meaningful.
Indeed, $\pi$~is \Lip\ continuous and, 
by combining \eqref{regphi} and \eqref{regsigma} with our assumptions on the structure and the data,
one can show that
(similarly as in the proof of the forthcoming estimates~\accorpa{stimagammakappa}{stimah})
\Bsist
  && \gamma(\phi) \, \sigma \,, \kappa(\phi) 
  \in \L\infty H ,
  \label{reggammakappa}
  \\[2mm]
  && h(\phi) \bigl( m(\sigma) - \mz u \bigr) \, \phi 
  \in \pier{\L2{\Lx\pz}}, 
  \label{reghphi}
  \\[2mm]
  && h(\phi) \bigl( m(\sigma) - \mz u \bigr)
  \in \pier{\L2{\Lx\ph}},
  \label{reghv}
\Esist
and we observe that every test function $v$ in \eqref{prima} belongs to $\Lx\php$ by \eqref{hpembedding}.
Finally, we stress that \eqref{prima} and \eqref{seconda} 
are equivalent to their time-integrated variants with time dependent test functions.
For instance, the former is equivalent~to
\Bsist
  && \intQ \dt\phi (\phi - v)
  + \ioT \bigl( A^\rho\phi(t) , A^\rho (\phi(t) - v(t)) \bigr) \, dt
  + \intQ \Beta(\phi)
  + \intQ \pi(\phi) ( \phi - v)
  \non
  \\
  && \leq \intQ h(\phi) \bigl( m(\sigma)) - \mz \, u \bigr) (\phi - v)
  + \intQ \Beta(v)
  \quad \hbox{for every $v\in\L2{\VA\rho}$},
  \qquad
  \label{intprima}
\Esist
\pier{with the same warning as above for the last term.}  
Also in this inequality and in the analogous equation for~$\sigma$, 
all the integral are meaningful
due to \Regsoluz\ and \accorpa{reggammakappa}{reghv}.
Here is our existence result.

\Bthm
\label{Existence}
Let the assumptions \HPstruttura\ on the structure of the system
and \HPdati\ on the data be fulfilled.
Then there exists at least \juerg{one} pair $(\phi,\sigma)$ satisfying \Regsoluz\ and solving  problem \Pbl. 
Moreover, such a solution can be constructed \juerg{that} satisfies the estimate
\Beq
  \norma\phi_{\H1H\cap\L\infty{\VA\rho}}
  + \norma\sigma_{\H1H\cap\L\infty{\VB\tau}}
  + \norma{\Beta(\phi)}_{\L\infty\Luno}
  \leq \overline C_1,
  \label{stimasoluz}
\Eeq
with a constant $\overline C_1$ that depends only 
on the structure of the system, the norms of the data corresponding to \HPdati, and~$T$.
\Ethm

In order to prove uniqueness and continuous dependence, we have to reinforce our assumptions on the structure.
Namely, we \juerg{make the following postulates:}
\Bsist
  && \phu \,,\, \qhu \,,\, \phd \,,\, \qhd \,,\, \pgu \,,\, \qgu \,,\, \rgu \,,\, \pgd \,,\, \qgd \,,\, \pk \,,\, \qk
  \in [2,+\infty)
  \quad \hbox{satisfy}
  \non
  \\[2mm]
  && \quad \frac 1\phu + \frac 1\qhu
  = \frac 1\phd + \frac 1\qhd
  = \frac 1\pgu + \frac 1\qgu + \frac 1\rgu
  \non
  \\
  && \quad {}
  = \frac 1\pgd + \frac 1\qgd
  = \frac 1\pk + \frac 1\qk
  = \frac 12\,.
  \label{hpparambis}
  \\
  && \VA\rho \subset \Lx\pstar
  \aand
  \VB\tau \subset \Lx\qstar
  \quad \hbox{with continuous embeddings, where} 
  \qquad
  \non
  \\
  && \quad \pstar := \max \{
    \phu(\ah-1) \,,\, \qhu \,,\, \ah\phd \,,\, \pgu(\ag-1) \,,\, \qgu \,,\, \ag\pgd \,,\, \pk \,,\, \qk
  \}
  \non
  \\
  && \aand
  \qstar := \max \{ \qhd, \rgu, \qgd \} \,.
  \label{hpembeddingbis}
  \\
  \separa
  && \hbox{$m$ is \Lip\ continuous\,.}
  \label{hpmbis}
  \\
  \separa
  && \hbox{$h$, $\gamma$ and $\kappa$ are of class $C^1$ and satisfy}
  \non
  \\
  && \quad |h'(r)| \leq C_0' \, |r|^{\ah-1} + C_1' \,,
  \quad |\gamma'(r)| \leq C_0' \, |r|^{\ag-1} + C_1'\,,
  \non
  \\
  && \aand |\kappa'(r)| \leq C_0' \, |r|^{\ak-1} + C_1'\,,
  \quad \hbox{for every $r\in\erre$,} 
  \label{hpgrowthbis}
\Esist
where $C_0'$ and $C_1'$ are given constants.
We notice that the inequalities \juerg{\eqref{hpgrowthbis}} imply both \eqref{hpgrowth} and the inequality
\Beq
  |h(r) - h(s)|
  \leq  \bigl( 
    C_0' \max\{ |r|^{\ah-1},|s|^{\ah-1} \} + C_1' 
  \bigr) |r-s|
  \quad \hbox{for every $r,s\in\erre$}\,,
  \qquad
  \label{dahpgrowthbis}
\Eeq
as well as \juerg{its} analogues for $\gamma$ and~$\kappa$.

\Brem
\label{Remhpbis}
The assumptions \accorpa{hpparambis}{hpembeddingbis} look very complicated.
However, they are satisfied in a number of situations.
One is that of the system \accorpa{citeprima}{citecauchy} described in the Introduction,
as one immediately sees by also accounting for Remark~\ref{Partcase}.
A~nontrivial case is given by Example~\ref{Example2}, where we recall that
$\VA\rho\subset\Lx5$, $\VB\sigma\subset\Lx6$ and $\ah=\ag=\ak=5/4$.
Indeed, an admissible choice of the new parameters is the following:
$\pgu=15/2$, $\qgu=5$, $\rgu=6$, and all of the other exponents are \juerg{taken as}~$4$.
\Erem

Here is our result.

\Bthm
\label{Contdep}
Besides the hypotheses of Theorem~\ref{Existence},
assume \juerg{that} \accorpa{hpparambis}{hpgrowthbis} \juerg{are satisfied}.
Then the solution to problem \Pbl\ is unique.
Moreover, given a constant~$M$, 
let $u_i$, $S_i$ and $\phi_{0,i}$, $i=1,2$, ,
be two choices of $u$, $S$ and~$\phiz$, respectively,
and $(\phi_i,\sigma_i)$ be corresponding solutions,
and assume that
\Beq
  \norma{u_i}_{\L2\Linfty} \,, \
  \norma{S_i}_{\L2\Linfty} \,, \
  \norma{\phi_i}_{\L\infty{\VA\rho}} \,, \ 
  \norma{\sigma_i}_{\L\infty{\VB\tau}}
  \leq M
  \label{upperbound}
\Eeq
for $i=1,2$.
Then the estimate
\Bsist
  && \norma{\phi_1-\phi_2}_{\L\infty H\cap\L2{\VA\rho}}
  + \norma{\sigma_1-\sigma_2}_{\L\infty H\cap\L2{\VB\tau}}
  \non
  \\
  && \leq \overline C_2 \bigl(
    \norma{u_1-u_2}_{\L2\Linfty}
    + \norma{S_1-S_2}_{\L2\Linfty}
    + \norma{\phi_{0,1}-\phi_{0,2}}
  \bigr) 
  \quad 
  \label{contdep}
\Esist
holds true with a constant $\overline C_2$ that only depends on the structure of our system,~$T$, and~$M$.
\Ethm

The remainder of the paper is organized as follows.
The uniqueness and continuous dependence result is proved in Section~\ref{UNIQUENESS},
while the existence of a solution is established
in the last Section~\ref{EXISTENCE} and is prepared by the study 
of the approximating problem introduced in Section~\ref{APPROX}.

\vspace{3mm}
Throughout the paper, we make a wide use of the \Holder\ inequality 
and of the elementary Young inequality
\Beq
  ab\leq \delta a^2 + \frac 1 {4\delta}\,b^2
  \quad \hbox{for every $a,b\in\erre$ and $\delta>0$} \,.
  \label{young}
\Eeq
Moreover, we use the notation
\Beq
  Q_t := \Omega \times (0,t)
  \quad \hbox{for $t\in(0,T)$},
  \label{defQt}
\Eeq
\pier{so that $Q=Q_T$.}
Finally, we state a general rule concerning the constants we are going to follow.
We always use a small-case italic $c$ without subscripts
for~different constants that may only depend the structure of our system 
(i.e.,~the operators $A^\rho$ and~$B^\tau$,  
the shape of the nonlinearities, the parameters that appear in our assumptions),
the final time~$T$
and the properties of the data involved in the statements at hand.
In particular, the values of such constants do not depend 
on the regularization parameter $\eps$ we introduce in Section~\ref{APPROX}.
Symbols like~$\cdelta$ (e.g., with $\delta=\eps$) denote constants that depend on the parameter~$\delta$, in addition.
It has to be clear that the values of $c$ and $\cdelta$ might change from line to line 
and even within the same formula or chain of inequalities. 
In contrast, we use different symbols (e.g., capital letters like $C_i$ in~\eqref{hpgrowth})
for precise values of constants we want to refer~to.


\section{Continuous dependence and uniqueness}
\label{UNIQUENESS}
\setcounter{equation}{0}

This section is devoted to the proof of the uniqueness and the continuous dependence
stated in Theorem~\ref{Contdep}.
More precisely, we prove just the continuous dependence,
since uniqueness follows as a consequence.
We pick two choices of the data as in the statement, the corresponding solutions
and a constant $M$ satisfying~\eqref{upperbound}.
We set for convenience
\Beq
  u := u_1-u_2 \,, \ \
  S := S_1-S_2 \,, \ \
  \phiz := \phi_{0,1}-\phi_{0,2} \,, \ \
  \phi := \phi_1-\phi_2 
  \ \ \hbox{and} \ \
  \sigma := \sigma_1-\sigma_2 \,.
  \label{diff}
\Eeq 
We also set for brevity
\Beq
  \phistar:= \max \{|\phi_1|,|\phi_2|\} \quad \hbox{\pier{pointwise}}
  \label{defphibar}
\Eeq
and denote by $D$ the maximum of the norms of the embeddings \eqref{hpembedding} and~\eqref{hpembeddingbis}.
We assume that $D\geq1$ without loss of generality.
At this point, we are ready to start.
According to our general rule, we use the same symbols $c$ and~$\cdelta$ 
(where $\delta>0$ is chosen later~on)
for (possibly) different constants,
as explained at the end of the previous section. 
In this proof, the values of $c$ and $\cdelta$ are allowed to depend on~$M$, in addition.
We write \eqref{prima} at the time $s$ for $(\phi_1,\sigma_1)$ and $(\phi_2,\sigma_2)$
and test the inequality obtained by $\phi_2(s)$ and~$\phi_1(s)$, respectively.
Then, we sum up, integrate over $(0,t)$ with respect to~$s$
and add the same quantity $\intQt|\phi|^2$ to both sides.
Since the terms involving $\Beta$ cancel each other, we obtain~that
\Bsist
  && \frac 12 \, \norma{\phi(t)}^2
  + \iot \norma{\phi(s)}_{A,\rho}^2 \, ds
  \non
  \\
  \separa
  && \leq \frac 12 \, \norma\phiz^2
  + \intQt |\phi|^2
  - \intQt \bigl( \pi(\phi_1)  - \pi(\phi_2) \bigr) \phi
  \non
  \\
  && \quad {}
  + \intQt \bigl\{
    h(\phi_1) \bigl( m(\sigma_1) - \mz u_1 \bigr) 
    - h(\phi_2) \bigl( m(\sigma_2) - \mz u_2 \bigr) 
  \bigr\} \phi \,.
  \label{diffprima}
\Esist
By the Young inequality and the \Lip\ continuity of~$\pi$, we immediately see~that
\Beq
  - \intQt \bigl( \pi(\phi_1) - \pi(\phi_2) \bigr) \phi
  \leq c \intQt |\phi|^2 \,.
  \label{diffpi}
\Eeq
For the next term, we owe to both the assumptions of Theorem~\ref{Existence}
and the supplementary assumptions on the structure (in particular, to~\eqref{dahpgrowthbis}).
We have~that
\Bsist
  && \intQt \bigl\{
    h(\phi_1) \bigl( m(\sigma_1) - \mz u_1 \bigr) 
    - h(\phi_2) \bigl( m(\sigma_2) - \mz u_2 \bigr) 
  \bigr\} \phi
  \non
  \\
  \separa
  && \leq \intQt |h(\phi_1) - h(\phi_2)| \, |m(\sigma_1) - \mz u_1| \, |\phi|
  \non
  \\
  && \quad {}
  + \intQt |h(\phi_2)| \, |m(\sigma_1) - m(\sigma_2) - \mz u| \, |\phi|
  \label{diffh}
\Esist
and we estimate the last two integrals, separately.
We have~that
\Bsist
  && \intQt |h(\phi_1) - h(\phi_2)| \, |m(\sigma_1) - \mz u_1| \, |\phi|
  \leq c \intQt (\phistar^{\ah-1} + 1) (1 + |u_1|) |\phi|^2
  \non
  \\
  \separa
  && \leq c \iot (\norma{(\phistar(s))^{\ah-1}}_\phu + 1) (1 + \norma{u_1(s)}_\infty) \norma{\phi(s)}_\qhu \, \norma{\phi(s)} \, ds
  \non
  \\
  && \leq \delta \iot (\norma{(\phistar(s))^{\ah-1}}_\phu^2 + 1) \norma{\phi(s)}_\qhu^2 \, ds
  + \cdelta \iot (1 + \norma{u_1(s)}_\infty^2) \norma{\phi(s)}^2 \, ds
  \non
  \\
  \separa
  && \leq \delta \, D^2 (\norma\phistar_{\L\infty{\Lx{\pier{\pstar}}}}^{2(\ah-1)} + 1) \iot \norma{\phi(s)}_{A,\rho}^2 \, ds
  + \cdelta \iot (1 + \norma{u_1(s)}_\infty^2) \norma{\phi(s)}^2 \, ds
  \non
  \\
  && \leq \delta \, D^2 \, D^{2(\ah-1)} (\norma\phistar_{\L\infty{\VA\rho}}^{2(\ah-1)} + 1)
     \iot \norma{\phi(s)}_{A,\rho}^2 \, ds
  \non
  \\
  && \quad {}
  + \cdelta \iot (1 + \norma{u_1(s)}_\infty^2) \norma{\phi(s)}^2 \, ds
  \non
  \\
  \separa
  && \leq \delta \, D^{2\ah} (M^{2(\ah-1)} + 1) \iot \norma{\phi(s)}_{A,\rho}^2 \, ds
  + \cdelta \iot (1 + \norma{u_1(s)}_\infty^2) \norma{\phi(s)}^2 \, ds \,.
  \label{diffh1}
\Esist
The other integral is estimated this~way
(as~for $\ph$, recall \accorpa{hpparam}{hpembedding}):
\Bsist
  && \intQt |h(\phi_2)| \, |m(\sigma_1) - m(\sigma_2) - \mz u| \, |\phi|
  \non
  \\
  && \leq c \intQt (|\phi_2|^\ah + 1) |\sigma| \, |\phi|
  + c  \intQt (|\phi_2|^\ah + 1) |u| \, |\phi|
  \non
  \\
  \separa
  && \leq c \iot (\norma{|\phi_2(s)|^\ah}_\phd + 1) \norma{\sigma(s)}_\qhd \, \norma{\phi(s)} \, ds
  \non
  \\
  && \quad {}
  + c  \iot (\norma{|\phi_2(s)|^\ah}_\ph + 1) \norma{u(s)}_\infty \, \norma{\phi(s)}_\php \, ds
  \non
  \\
  \separa
  && \leq \delta \iot \norma{\sigma(s)}_\qhd^2 \, ds
  + \cdelta (\norma{|\phi_2|^\ah}_{\juerg{L^\infty(0,T;L^{p_{h,2}}(\Omega))}}^2 + 1) \iot \norma{\phi(s)}^2 \, ds
  \non
  \\
  && \quad {}
  + \delta \iot \norma{\phi(s)}_\php^2 \, ds
  + \cdelta \, (\norma{|\phi_2|^\ah}_{\L\infty{\Lx\ph}}^2 + 1) \iot \norma{u(s)}_\infty^2 \, ds
  \non
  \\
  \separa
  && \leq \delta \, D^2 \iot \norma{\sigma(s)}_{B,\tau}^2 \, ds
  + \cdelta (\norma{\phi_2}_{\L\infty{\Lx{\ah\phd}}}^{2\ah} + 1) \iot \norma{\phi(s)}^2 \, ds
  \non
  \\
  && \quad {}
  + \delta \, D^2 \iot \norma{\phi(s)}_{A,\rho}^2 \, ds
  + \cdelta (\norma{\phi_2}_{\L\infty{\Lx{\ah\ph}}}^{2\ah} + 1) \, \norma u_{\L2\Linfty}^2 
  \\
  \separa
  && \leq \delta \, D^2\iot \norma{\sigma(s)}_{B,\tau}^2 \, ds
  + \cdelta \iot \norma{\phi(s)}^2 \, ds
  \non
  \\
  && \quad {}
  + \delta \, D^2 \iot \norma{\phi(s)}_{A,\rho}^2 \, ds
  + \cdelta \, \norma u_{\L2\Linfty}^2 \,.
  \label{diffh2}
\Esist
By collecting \accorpa{diffprima}{diffh2},
we conclude~that
\Bsist
  && \frac 12 \, \norma{\phi(t)}^2
  + \iot \norma{\phi(s)}_{A,\rho}^2 \, ds
  \non
  \\
  \separa
  && \leq \frac 12 \, \norma\phiz^2
  + \delta \bigl\{
    D^{2\ah} (M^{2(\ah-1)} + 1)
    + D^2
  \bigr\} \iot \norma{\phi(s)}_{A,\rho}^2 \, ds
  \non
  \\
  && \quad {}  
  + \cdelta \, \norma u_{\L2\Linfty}^2
  + \cdelta \iot (1 + \norma{u_1(s)}_\infty^2) \norma{\phi(s)}^2 \, ds
	\non
	\\
	&&\juerg{\quad{}
	+ \delta D^2\iot \|\sigma(s)\|^2_{B,\tau}\,ds\,,}
  \label{conclprima}
\Esist
and we observe that the function $s\mapsto\norma{u_1(s)}_\infty^2$ belongs to $L^1(0,T)$ 
and that its norm is bounded by~$M^2$.
Now, we write \eqref{seconda} at the time $s$ for both solutions,
test the difference by $\sigma(s)$, integrate over $(0,t)$ 
and add the same quantity $\intQt|\sigma|^2$ to both sides.
We obtain~that
\Bsist
  && \frac 12 \, \norma{\sigma(t)}^2 
  + \iot \norma{\sigma(s)}_{B,\tau}^2 \, ds
  \non
  \\
  && = \intQt |\sigma|^2
  - \intQt \bigl( \gamma(\phi_1) \sigma_1 - \gamma(\phi_2) \sigma_2 \bigr) \sigma
  \non
  \\
  && \quad {}
  + \intQt \bigl( \kappa(\phi_1) - \kappa(\phi_2) \bigr) \sigma
  - \intQt (S_1 \phi_1 - S_2 \phi_2) \sigma \,.
  \label{diffseconda}
\Esist
We estimate the last three terms, separately,
by accounting for the analogues of~\eqref{dahpgrowthbis} for $\gamma$ and~$\kappa$.
As for the first one,
we have that
\Bsist
  && - \intQt \bigl( \gamma(\phi_1) \sigma_1 - \gamma(\phi_2) \sigma_2 \bigr) \sigma
   = - \intQt \bigl( 
    (\gamma(\phi_1) - \gamma(\phi_2)) \sigma_1 
    + \gamma(\phi_2) \sigma
   \bigr) \sigma
  \non
  \\
  \separa
  && \leq c \intQt (\phistar^{\ag-1} + 1) |\phi| \, |\sigma_1| \, |\sigma|
  + c \intQt (|\phi_2|^\ag + 1) \, |\sigma|^2
  \non
  \\
  \separa
  && \leq c \iot (\norma{(\phistar(s))^{\ag-1}}_\pgu + 1) \, \norma{\phi(s)}_\qgu \, \norma{\sigma_1(s)}_\rgu \, \norma{\sigma(s)} \, ds
  \non
  \\
  && \quad {}
  + c \iot (\norma{(\phi_2(s))^\ag}_\pgd + 1) \, \norma{\sigma(s)}_\qgd \, \norma{\sigma(s)} \, ds
  \non
  \\
  \separa
  && \leq \delta \iot (\norma{\phistar(s)}_{\pier{\pstar}}^{2(\ag-1)} + 1) \, \norma{\phi(s)}_\qgu^2 \, ds
  + \cdelta \iot \norma{\sigma_1(s)}_\rgu^2 \, \norma{\sigma(s)}^2 \, ds
  \non
  \\
  && \quad {}
  + \delta \iot \norma{\sigma(s)}_\qgd^2 \, ds
  + \cdelta \iot (\norma{\phi_2(s)}_{\pgd\ag}^{2\ag} + 1) \, \norma{\sigma(s)}^2 \, ds
  \non
  \\
  \separa
  && \leq \delta \, D^{2\ag} \iot (\norma{\phistar(s)}_{A,\rho}^{2(\ag-1)} + 1) \, \norma{\phi(s)}_{A,\rho}^2 \, ds
  + \cdelta \, D^2 \iot \norma{\sigma_1(s)}_{B,\tau}^2 \, \norma{\sigma(s)}^2 \, ds
  \non
  \\
  && \quad {}
  + \delta \, D^2 \iot \norma{\sigma(s)}_{B,\tau}^2 \, ds
  + \cdelta \, D^{2\ag} \iot (\norma{\phi_2(s)}_{A,\rho}^{2\ag} + 1) \, \norma{\sigma(s)}^2 \, ds
  \non
  \\
  \separa
  && \leq \delta \, D^{2\ag} \, (M^{2(\ag-1)} + 1) \iot \norma{\phi(s)}_{A,\rho}^2 \, ds
  + \cdelta \iot \norma{\sigma(s)}^2 \, ds
  \non
  \\
  && \quad {}
  + \delta \, D^2 \iot \norma{\sigma(s)}_{B,\tau}^2 \, ds
  + \cdelta \iot \norma{\sigma(s)}^2 \, ds \,.
  \label{diffgamma}
\Esist
We estimate the next term in this way:
\Bsist
  && \intQt \bigl( \kappa(\phi_1) - \kappa(\phi_2) \bigr) \sigma
  \leq c \intQt (\phistar^{\ak-1} + 1) \, |\phi| \, |\sigma|
  \non
  \\
  && \leq c \iot (\norma{(\phistar(s))^{\ak-1}}_\pk + 1) \, \norma{\phi(s)}_\qk \, \norma{\sigma(s)} \, ds
  \non
  \\
  \separa
  && \leq \delta \iot (\norma{(\phistar(s))^{\ak-1}}_\pk^2 + 1) \, \norma{\phi(s)}_\qk^2 \, ds
  + \cdelta \iot \norma{\sigma(s)}^2 \, ds
  \non
  \\
  \separa
  && \leq \delta \, D^{2\ak} \iot (\norma{\phistar(s)}_{A,\rho}^{2(\ak-1)} + 1) \, \norma{\phi(s)}_{A,\rho}^2 \, ds
  + \cdelta \iot \norma{\sigma(s)}^2 \, ds
  \non
  \\
  \separa
  && \leq \delta \, D^{2\ak} \, (M^{2(\ak-1)} + 1) \iot\norma{\phi(s)}_{A,\rho}^2 \, ds
  + \cdelta \iot \norma{\sigma(s)}^2 \, ds \,.
  \label{diffkappa}
\Esist
Finally, \pier{we observe that}
\Bsist
  && - \intQt (S_1 \phi_1 - S_2 \phi_2) \sigma
  = - \intQt \bigl(
    S \, \phi_1
    + S_2 \, \phi
  \bigr) \sigma
  \non
  \\
  \separa
  && \leq \intQt |S| \, |\phi_1| \, |\sigma|
  + \intQt |S_2| \, |\phi| \, |\sigma|
  \non
  \\
  && \leq c \iot \norma{S(s)}_\infty \, \norma{\phi_1(s)} \, \norma{\sigma(s)} \, ds
  +  c \iot \norma{S_2(s)}_\infty \, \norma{\phi(s)} \, \norma{\sigma(s)} \, ds
  \non
  \\
  \separa
  && \leq c \ioT \norma{S(s)}_\infty^2 \, ds
  + c \iot \norma{\sigma(s)}^2 \, ds
  \non
  \\
  && \quad {}
  + c \iot \norma{\phi(s)}^2 \, ds
  + c \iot \norma{S_2(s)}_\infty^2 \,\norma{\sigma(s)}^2 \, ds \,.
  \label{diffS}
\Esist
By collecting \accorpa{diffseconda}{diffS} and rearranging,
we conclude~that
\Bsist
  && \frac 12 \, \norma{\sigma(t)}^2 
  + \iot \norma{\sigma(s)}_{B,\tau}^2 \, ds
  \non
  \\
  && \leq \delta \, \bigl\{
    D^{2\ag} \, (M^{2(\ag-1)} + 1) 
    + D^{2\ak} \, (M^{2(\ak-1)} + 1)
  \bigr\} \iot \norma{\phi(s)}_{A,\rho}^2 \, ds
  \non
  \\
  && \quad {}
  + \delta \, D^2 \iot \norma{\sigma(s)}_{B,\tau}^2 \, ds
  + c \, \norma S_{\L2\Linfty}^2 
  \non
  \\
  && \quad {}
  + c \iot \norma{\phi(s)}^2 \, ds
  + \cdelta \iot (1 + \norma{S_2(s)}_\infty^2) \,\norma{\sigma(s)}^2 \, ds\,,
  \label{conclseconda}
\Esist
where we observe that the function $s\mapsto\norma{S_2(s)}_\infty^2$ belongs to $L^1(0,T)$
and that its norm is bounded by~$M^2$. 
At this point, we add \eqref{conclprima} and \eqref{conclseconda} to each other,
choose $\delta$ small enough and apply the Gronwall lemma.
This yields \eqref{contdep} with a constant that has the same dependence as the constant $\overline C_2$ in the statement.
This completes the proof.


\section{Approximation}
\label{APPROX}
\setcounter{equation}{0}

In this section, we deal with an approximation of problem \Pbl\ depending on the parameter $\eps\in(0,1)$
and solve it by a Faedo--Galerkin scheme.
We first introduce the Moreau--Yosida regularizations 
$\Betaeps$ and $\betaeps$ of $\Beta$ \juerg{and} $\beta$ at the level~$\eps>0$
(see, e.g., \cite[p.~28 and p.~39]{Brezis}),
and we recall that $\betaeps$ is the derivative of $\Betaeps$ and \pier{is} \Lip\ continuous.
We also remark that the following properties hold true (the first inequality being due to~\eqref{hpBeta}):
\begin{align}
  & 0 \leq \Beta_{\eps''}(r) \leq \Beta_{\eps'}(r) \leq \Beta(r)
  \quad \hbox{if $0<\eps'\leq\eps''$}
  \non
  \\
  & \aand
  \lim_{\eps\searrow0} \Betaeps(r) = \Beta(r),
  \quad \hbox{for every $r\in\erre$}\,.
  \label{propBetaeps}
\end{align}
The approximating problem, whose unknown is the pair $(\phieps,\sigmaeps)$,
is obtained by replacing $\Beta$ by~$\Betaeps$ in \Pbl\
and approximating the other nonlinearities by smoother functions.
Moreover, we also regularize the data $u$ and~$S$.
However, since $\Betaeps$ is differentiable, the variational inequality (i.e., the analogue of \eqref{prima})
can be replaced by an equation involving the derivative~$\betaeps$ of~$\Betaeps$.
As for the other approximating nonlinearities, we require that they are bounded and \Lip\ continuous 
and converge to the original ones uniformly on every compact interval of~$\erre$.
Moreover, inequalities that are analogous to \accorpa{hpgrowth}{hphBeta} 
should be satisfied by the regularized functions uniformly with respect to~$\eps$.
A~possible construction of such approximations is based on the lemma stated below.
It is clear that \juerg{the second assumption in} \eqref{hpKi} is empty
if \juerg{$D(\Beta)=\erre$}.
In this case, $K_3$~can be any positive constant.

\Blem
\label{Approx}
Let $\psi:\erre\to\erre$ be continuous and define $\psi^\eps,\psieps:\erre\to\erre$ by the formulas
\begin{align}
  & \psi^\eps(r) := \psi(r) \quad \hbox{if $|r|\leq1/\eps$} \,,
  \quad \psi^\eps(r) := \psi((\sign r)/\eps) \quad \hbox{if $|r|>1/\eps$}
  \non
  \\
  & \aand
  \psieps(r) := \frac 1{2\eps} \int_{r-\eps}^{r+\eps} \psi^\eps(s) \, ds
  \quad \hbox{for $r\in\erre$}.
  \label{defpsieps}
\end{align}
Then, $\psieps$ is bounded and \Lip\ continuous \juerg{on $\erre$,
and it holds} that
\Beq
  \psieps \to \psi
  \quad \hbox{uniformly on every bounded interval of $\erre$}.
  \label{convpsieps}
\Eeq
Moreover, if 
\Beq
  |\psi(r)| \leq K_0 \, |r|^\alpha + K_1
  \quad \hbox{for every $r\in\erre$ and some positive constants $\alpha$, $K_0$ and $K_1$}\,,
  \label{disugpsi}
\Eeq
then there are constants $\hat K_0$ and $\hat K_1$ depending only on $\alpha$, $K_0$ and $K_1$ such~that 
\Beq
  |\psieps(r)| \leq \hat K_0 \, |r|^\alpha + \hat K_1
  \quad \hbox{for every $r\in\erre$ and $\eps\in(0,1)$} \,.
  \label{disugpsieps}
\Eeq
Finally, also assume that
\Beq
  |\psi(r)|^2 \leq K_2 \, \Beta(r) + K_3
  \quad \hbox{for every $r\in\erre$},
  \label{disugpsiBeta}
\Eeq
where $K_2$ and $K_3$ are constants satisfying
\Beq
  K_2 \geq 1
  \aand
  K_3 \geq \sup_{|r-\rz|\leq\delta} |\psi(r)|^2
  \label{hpKi}
\Eeq
for every end-point $\rz$ of $D(\Beta)$ (if any) and some $\delta>0$.
Then, the function $\psitilde:\erre\to\erre$ defined by setting
\Beq
  \psitilde(r) := 
  \max \bigl\{
    - (K_2\,\Betaeps(r)+K_3)^{1/2} ,
    \min \{\psieps(r) , (K_2\,\Betaeps(r)+K_3)^{1/2} \}
  \bigr\} 
  \label{defpsitilde}
\Eeq
is bounded and \Lip\ continuous \juerg{on $\erre$} and satisfies
\Beq
  |\psitilde(r)|^2 \leq K_2 \, \Betaeps(r) + K_3
  \quad \hbox{for every $r\in\erre$ and $\eps\in(0,1)$},
  \label{disugpsitilde} 
\Eeq
as well as the analogues of \eqref{convpsieps} and~\eqref{disugpsieps}.
\Elem

\Bdim
Clearly, $\psieps$ is of class $C^1$ and its derivative is given by $(1/(2\eps)(\juerg{\psi}^\eps(r+\eps)-\psi^\eps(r-\eps)$.
In particular, if $|r|>(1/\eps)+\eps$, we have that $\psieps'(r)=0$, 
so that $\psieps$ is bounded and \Lip\ continuous \juerg{on $\erre$}.
Take now any $M>0$ and $\eta>0$.
Since $\psi$ is uniformly continuous on every bounded interval, 
we have that
$|\psi(r)-\psi(s)|\leq\eta$ whenever 
$\eps\in(0,1)$ \juerg{is} small enough, $|r|,|s|\leq M+1\,$ and $\,|r-s|\leq\eps$.
Then, if $\eps$ also satisfies $M+1<1/\eps$ 
and $r$ belongs to $[-M,M]$,
we conclude that
\Beq
  |\psieps(r) - \psi(r)|
  = \Bigl| \frac 1{2\eps} \int_{r-\eps}^{r+\eps} \bigl( \psi(s) - \psi(r) \bigr) \, ds \Bigr|
  \leq \frac 1{2\eps} \int_{r-\eps}^{r+\eps} \eta \, ds
  = \eta \,.
  \non
\Eeq
All this proves \eqref{convpsieps}.
Assume now \eqref{disugpsi}.
Then, if $|s|\leq|r|$, we have that 
$|\psi(s)|\leq K_0\,|s|^\alpha+K_1\leq K_0\,|r|^\alpha+K_1$.
Therefore, 
$\sup_{|s|\leq|r|}|\psi(s)|\leq K_0\,|r|^\alpha+K_1$,
whence also
\Beq
  \sup_{|s|\leq|r|} |\psi^\eps(s)|
  \leq \sup_{|s|\leq|r|} |\psi(s)|
  \leq K_0 \, |r|^\alpha + K_1 \,.
  \non
\Eeq
Set now $M_\alpha:=\sup_{r\in\erre}(|r|+1)^\alpha/\juerg{(|r|^\alpha+1)}$\juerg{, which is obviously finite}.
Then, for every $r\in\erre$ and $\eps\in(0,1)$, we have that
\Beq
  |\psieps(r)|
  \leq \sup_{|s|\leq|r|+1} |\psi^\eps(s)|
  \leq K_0 (|r|+1)^\alpha + K_1
  \leq K_0 M_\alpha (|r|^\alpha + 1) + K_1\,,
  \non
\Eeq
so that \eqref{disugpsieps} holds with $\hat K_0=K_0M_\alpha$ and $\hat K_1=K_0M_\alpha+K_1$.

Let us come to the properties of $\psitilde$ under the assumption~\eqref{disugpsiBeta}.
First of all, recall that $K_2$ and $K_3$ are positive, that $\Betaeps$~is nonnegative and 
locally \Lip\ continuous, and that it tends to $+\infty$ as its argument tends to~$\pm\infty$
(as~a consequence of the last \pier{condition in}~\eqref{hpBeta}).
Thus, \pier{it turns out that} $\psitilde(r)=\psieps(r)$ for $|r|$ large enough. 
Moreover, $\psieps$~is (globally) \Lip\ continuous. 
This yields that the function $\psitilde$ it is well defined and \Lip\ continuous.
Moreover, $\psitilde$~is bounded \juerg{in view of the boundedness of $\psieps$}.
Furthermore, \eqref{disugpsitilde} trivially follows \juerg{from} the definition of~$\psitilde$,
and the analogue of \eqref{disugpsieps} for $\psitilde$ is a consequence of \eqref{disugpsieps} itself,
since $|\psitilde|\leq|\psieps|$.
Let us prove the analogue of~\eqref{convpsieps}
by assuming that the interior of $D(\Beta)$ is nonempty
(the opposite case is even easier since \juerg{then} $D(\Beta)$~is a singleton).
Take a compact interval $I$ contained in the interior of~$D(\Beta)$.
Then, the restriction of $\Beta$ to $I$ is continuous.
Therefore, by \eqref{propBetaeps} and Dini's theorem on monotone convergence,
$\Betaeps$~converges to $\Beta$ uniformly in~$I$.
By combining this with \eqref{convpsieps}, 
we infer that $\psitilde$ uniformly converges in $I$ to the function 
\Beq
  I \ni r \mapsto
  \max \bigl\{
    - (K_2\,\Beta(r)+K_3)^{1/2} ,
    \min \{\psi(r) , (K_2\,\Beta(r)+K_3)^{1/2} \}
  \bigr\} 
  = \psi(r) \,,
  \non
\Eeq
the last equality being due to \eqref{disugpsiBeta}.
If $D(\Beta)=\erre$, \juerg{then} the convergence properties under investigation \juerg{are} completely proved.
In the opposite case, we consider two situations regarding the compact interval~$I$.
In the first one, $I=[a,b]$~is contained in the exterior of~$D(\Beta)$.
Then, $\Betaeps(r)$~tends to $+\infty$ as $\eps$ tends to zero for every $r\in I$.
Moreover, the convergence is monotone due to \eqref{propBetaeps}.
By applying Dini's theorem to $1/\Betaeps$ (whose pointwise limit is obviously continuous),
we deduce that $\Betaeps$ is uniformly divergent.
Thus, by setting $M:=\sup_{a-1\leq r\leq b+1}|\psi(r)|$,
we have that $\Betaeps(r)\geq M^2$ for every $r\in I$ and $\eps>0$ small enough.
On the other hand, $|\psieps(r)|\leq M$ for every $r\in I$ (see the first part of this proof).
We thus have for such values of~$\eps$ 
(since $K_2\geq1$ and $\Betaeps$ is nonnegative) that
\Beq
  |\psieps(r)|^2
  \leq M^2
  \leq \Betaeps(r) 
  \leq K_2 \, \Betaeps(r) + K_3
  \quad \hbox{for every $r\in I$}.
  \non
\Eeq
Therefore, $\psitilde$ coincides with $\psieps$ on $I$ for $\eps>0$ small enough,
and our assertion follows from \eqref{convpsieps}.
In the other situation, $I=[\rz-\delta,\rz+\delta]$ is the $\delta$-\nbh\ of an \juerg{endpoint} $r_0$ of $D(\beta)$
like in~\eqref{hpKi}.
Then, for every $r\in I$, we have that $|\psieps(r)|^2\leq K_3\leq K_2\,\Betaeps(r)+K_3$,
whence also $\psitilde(r)=\psieps(r)$, so that
our assertion follows from~\eqref{convpsieps} also in this case.
Since every compact interval of $\erre$ is the union of $n\leq3$ intervals of the previous type,
the uniform convergence property we have claimed is completely proved.
\Edim


We use the above lemma to introduce the approximating nonlinearities.
We choose $K_0=C_0$ and $K_1=C_1$, the constants appearing in~\eqref{hpgrowth},
and let $\alpha$ take the values $\ah$, $\ag$, and $\ak$, according to the functions we want to define.
Finally, by recalling~\eqref{hphBeta}, we set $K_2:=\max\{C_2,1\}$ and, if $D(\Beta)\not=\erre$,
we choose $K_3\geq C_3$ in order to satisfy~\eqref{hpKi} with $\psi=h$
(in~particular, \eqref{disugpsiBeta} holds with $\psi=h$ as a consequence of~\eqref{hphBeta}).
Then, we \pier{agree~that}
\Bsist
  && \hbox{$\meps$,$\gammaeps$, and $\kappaeps$, are defined as $\psieps$ with $\psi=m,\gamma,\kappa$, respectively,}
  \non
  \\
  && \quad \hbox{and $\heps$ is defined as $\psitilde$ with $\psi=h$} \,.
  \label{defnonlineps}
\Esist
Finally, we replace the data $u$ and~$S$
by approximating data $\ueps$ and $\Seps$ satisfying
\Bsist
  && \ueps \in \LQ\infty
  \aand
  \Seps \in \LQ\infty
  \quad \hbox{with}
  \non
  \\
  && |\ueps| \leq |u| \,, \quad
  |\Seps| \leq |S| \,, \quad
  \ueps \to u
  \aand 
  \Seps \to S
  \quad \aeQ \,.
  \label{hpuSeps}
\Esist
To fulfill these conditions, one can set, e.g., 
$\ueps:=\max\{-1/\eps,\min\{u,1/\eps\}\}$, and analogously define~$\Seps$.
Therefore, the problem we consider is the following:
\Bsist
  && \iO \dt\phieps(t) \, v
  + \bigl( A^\rho\phieps(t) , A^\rho v \bigr)
  + \iO (\betaeps+\pi)(\phieps(t)) \, v
  \non
  \\
  && = \iO \heps(\phieps(t)) \bigl( \meps(\sigmaeps(t)) - \mz \, \ueps(t) \bigr) \, v
  \non
  \\
  && \quad \hbox{for every $v\in\VA\rho$ and \aat},
  \label{primaeps}
  \\[2mm]
  \separa
  && \iO \dt\sigmaeps(t) \, v
  + \bigl( B^\tau \sigmaeps(t) , B^\tau v \bigr)
  + \iO \gammaeps(\phieps(t)) \, \sigmaeps(t) \, v
  \non
  \\
  && = \iO \kappaeps(\phieps(t)) \, v
  - \iO \Seps(t) \phieps(t) \, v
  \non
  \\
  && \quad \hbox{for every $v\in\VB\tau$ and \aat},
  \label{secondaeps}
  \\
  && \phieps(0) = \phiz 
  \aand
  \sigmaeps(0) = \sigmaz \,.
  \label{cauchyeps}
\Esist
\Accorpa\Pbleps primaeps cauchyeps

\Brem
\label{RemPbleps}
We stress once more that \eqref{primaeps} is equivalent to the variational inequalities
obtained by replacing the nonlinearities with their approximations defined above in both \eqref{prima} and \eqref{intprima}.
Moreover, \eqref{primaeps} and \eqref{secondaeps} are also equivalent to 
their time-integrated versions with test functions taken in $\L2{\VA\rho}$ and $\L2{\VB\tau}$, respectively.
\Erem

\Bthm
\label{Wellposednesseps}
Under \pier{the same assumptions as in} Theorem~\ref{Existence},
the approximating problem \Pbleps\ has a unique solution $(\phieps,\sigmaeps)$
satisfying the analogues of the regularity conditions \accorpa{regphi}{regsigma}.
\Ethm

The remainder of the section is devoted to the proof of this theorem.
As for uniqueness, we can apply Theorem~\ref{Contdep}, 
since $\betaeps$ has the same properties \juerg{as} $\beta$,
and the functions $\meps$, $\heps$, $\gammaeps$ and $\kappaeps$ are bounded and have bounded derivatives,
so that all of the assumptions that are needed are satisfied by the approximating nonlinearities.
Hence, we just have to prove the existence of a solution.
To this end, we introduce a discrete problem depending on the parameter $n\in\enne$ by means of a Faedo--Galerkin scheme.
Then, we solve it and then take the limit of its solution as $n$ tends to infinity.

\step
The discrete problem

We recall that $e_j$ and $e'_j$, $j=1,2,\dots$, are the eigenfunctions of the operators $A$ and~$B$, respectively.
For every integer $n\geq1$ we~set
\Beq
  \VAn := \Span\graffe{e_1,\dots,e_n}
  \aand
  \juerg{\VBn} := \Span\graffe{e'_1,\dots,e'_n}
  \label{defspazin}
\Eeq
and look for a pair $(\phin,\sigman)$ enjoying the regularity
\Beq
  \phin \in \H1\VAn
  \aand
  \sigman \in H^1(0,T;\juerg{\VBn})
  \label{regsoluzn}
\Eeq
and solving the following problem
\Bsist
  && \iO \dt\phin(t) \, v
  + \bigl( A^\rho\phin(t) , A^\rho v \bigr)
  + \iO (\betaeps+\pi)(\phin(t)) \, v
  \non
  \\
  && = \iO \heps(\phin(t)) \bigl( \meps(\sigman(t)) - \mz \, \ueps(t) \bigr) \, v
  \non
  \\
  && \quad \hbox{for every $v\in\VAn$ and \aat},
  \label{priman}
  \\[2mm]
  \separa
  && \iO \dt\sigman(t) \, v
  + \bigl( B^\tau \sigman(t) , B^\tau v \bigr)
  + \iO \gammaeps(\phin(t)) \, \sigman(t) \, v
  \non
  \\
  && = \iO \kappaeps(\phin(t)) \, v
  - \iO \Seps(t) \phin(t) \, v
  \non
  \\
  && \quad \hbox{for every $v\in \juerg{\VBn}$ and \aat},
  \label{secondan}
  \\[2mm]
  && (\phin(0),v) = (\phiz,v)
  \aand
  (\sigman(0),v) = (\sigmaz,v)
  \non
  \\
  && \quad \hbox{for every $v\in\VAn$ and $v\in \juerg{\VBn}$, respectively}.
  \label{cauchyn}
\Esist
\Accorpa\Pbleps primaeps cauchyeps
Even though $\phin$ and $\sigman$ obviously depend on $\eps$ as well, we do not stress this in the notation.
We observe that \eqref{cauchyn} simply means that
\Beq
  \phin(0) = \somma j1n (\phiz,e_j) e_j
  \aand
  \sigman(0) = \somma j1n (\sigmaz,e'_j) e'_j \,,
  \label{cauchynbis}
\Eeq
since $\phin(0)\in\VAn$ and $\sigman(0)\in \juerg{\VBn}$.
Notice that this implies that
\Bsist
  && \norma{\phin(0)} \leq \norma\phiz \,, \quad
  \norma{\sigman(0)} \leq \norma\sigmaz \,, \quad
  \non
  \\
  && \norma{\phin(0)}_{A,\rho} \leq \norma\phiz_{A,\rho}
  \aand
  \norma{\sigman(0)}_{B,\tau} \leq \norma\sigmaz_{B,\tau} \,. 
  \label{dacauchynbis}
\Esist
Indeed, we have for instance that
\Bsist
  && \norma{A^\rho\phi(0)}^2
  = \Norma{{\textstyle\somma i1n} (\phiz,e_j) A^\rho e_j}^2
  = \Norma{{\textstyle\somma i1n} (\phiz,e_j) \lambda_j^\rho e_j}^2
  \non
  \\ 
  && \leq \Norma{{\textstyle\somma i1\infty} (\phiz,e_j) \lambda_j^\rho e_j}^2
  = \Norma{{\textstyle\somma i1\infty} (\phiz,e_j) A^\rho e_j}^2
  = \norma{A^\rho\phiz}^2 \,.
  \non
\Esist
The discrete problem has a unique solution, as we see at once.
By \eqref{regsoluzn}, the unknowns have to be expanded~as
\Beq
  \phin(t) = \somma j1n y_j(t) e_j
  \aand
  \sigman(t) = \somma j1n z_j(t) e'_j\,,
  \non
\Eeq
with some coefficients $y_j,z_j\in H^1(0,T)$.
Therefore, the true unknowns are the vectors $y:=(y_1,\dots,y_n)$ and $z:=(z_1,\dots,z_n)$.
Since it is sufficient to take $v=e_i$ and $v=e'_i$ 
with $i=1,\dots,n$ in \eqref{priman} and \eqref{secondan}, respectively,
and since the eigenvectors satisfy~\eqref{eigen},
the variational equations \eqref{priman} and \eqref{secondan} in terms of $y$ and $z$ become
\Bsist
  && y'_i(t) 
  + \lambda_i^{2\rho} \, y_i(t)
  + \Psi_{1,i}(t,y(t),z(t))
  = \Psi_{2,i}(t,y(t),z(t))
  \non
  \\
  && \quad \hbox{for $i=1,\dots,n$ and \aat},
  \label{primanbis}
  \\[2mm]
  && z'_i(t) 
  + (\lambda'_i)^{\pier{2\tau}} \, z_i(t)
  + \Psi_{3,i}(t,y(t),z(t))
  = \Psi_{4,i}(t,y(t),z(t))   
  \non
  \\
  && \quad \hbox{for $i=1,\dots,n$ and \aat},
  \label{secondanbis}
\Esist
where the Carath\'eodory functions $\Psi_{k,i}:(0,T)\times\erren\times\erren\to\erre$ are defined by
\Bsist
  && \Psi_{1,i}(t,\bar y,\bar z)
  := \iO (\betaeps+\pi) \Bigl( {\textstyle\somma j1n}\bar y_j e_j \Bigr) e_i\,,
  \non
  \\
  &&  \Psi_{2,i}(t,\bar y,\bar z)
  := \iO \heps \Bigl( {\textstyle\somma j1n} \bar y_j e_j \Bigr) 
    \Bigl[ \meps \Bigl( {\textstyle\somma j1n}\bar z_j e'_j \Bigr) - \mz \, \ueps(t) \Bigr] \, e_i\,,
  \non
  \\
  && \Psi_{3,i}(t,\bar y,\bar z)
  := \iO \gammaeps \Bigl( {\textstyle\somma j1n} \bar y_j e_j \Bigr) \, {\textstyle\somma j1n} \bar z_j e'_j \, e'_i\,,
  \non
  \\
  && \Psi_{4,i}(t,\bar y,\bar z)
  := \iO \kappaeps \Bigl( {\textstyle\somma j1n} \bar y_j e_j \Bigr) e'_i
  - \iO \Seps(t) {\textstyle\somma j1n} \bar y_j e_j \, e'_i \,,
  \non
\Esist
for $\bar y,\bar z\in\erren$ and \aat.
Hence, we have \pier{obtained} a system of $2n$ ordinary differential equations in $2n$ unknowns.
Since all of the functions $\betaeps,\dots,\kappaeps$ are \Lip\ continuous
and $\heps$, $\meps$ and $\gammaeps$ are even bounded, as well as $\ueps$ and~$\Seps$,
the functions $\Psi_{k,i}$ are \Lip\ continuous with respect to $(\bar y,\bar z)\in\erre^{2n}$
uniformly with respect to~$t$.
Since \eqref{cauchyn} (or \eqref{cauchynbis}) provides an initial condition for~$(y,z)$,
we conclude that the discrete problem has a unique solution with the regularity specified by~\eqref{regsoluzn}.

\medskip

\juerg{Now that} the discrete problem is solved, we can start estimating.
According our general rule, the symbol $\ceps$ 
denotes (possibly different) constants that are allowed to depend 
on the structure, the data, $T$, and~$\eps$, but not on~$n$.

\step
First a priori estimate

We write \eqref{priman} at the time $s$ and choose $v=\dt\phin(s)\in\VAn$.
Then, we integrate with respect to $s$ over the interval $(0,t)$ with $t\in(0,T]$.
Moreover, we add the same quantity 
$\frac12\norma{\phin(t)}^2-\frac12\norma{\phin(0)}^2=\intQt\phin\dt\phin$ 
to both sides.
After rearranging \juerg{terms}, we obtain \juerg{the identity}
\begin{align}
& \intQt |\dt\phin|^2 
  + \frac 12 \, \norma{\phin(t)}_{A,\rho}^2
  + \iO \Betaeps(\phin(t))
  \non
  \\
  & = \frac 12 \, \norma{\phin(0)}_{A,\rho}^2
  + \iO \Betaeps(\phin(0))
  \non
  \\
  &\quad{}+ \intQt \bigl(
    \phin - \pi(\phin)
    + \heps(\phin) \bigl( \meps(\sigman) - \mz \ueps \bigr) \bigr) \dt\phin \,.
  \label{pier1}
\end{align}
All of the terms on the \lhs\ are nonnegative,
and the whole \rhs\ can be estimated by the quantity
\Beq
  \ceps \, \norma{\phin(0)}_{A,\rho}^2
  + \frac 12 \iot |\dt\phin|^2
  + \ceps \intQt |\phin|^2
  + \ceps \,,
  \non
\Eeq
since $\Betaeps$ grows quadratically at infinity, $\pi$ is \Lip\ continuous \juerg{on $\erre$},
and the functions $\heps$, $\meps$, and $\ueps$, are bounded.
Moreover, we can owe to \eqref{dacauchynbis}.
We thus obtain the estimate
\Beq
  \norma\phin_{\H1H\cap\L\infty{\VA\rho}} \leq \ceps \,.
  \label{primastiman}
\Eeq

\step
Second a priori estimate

Similarly, we test \eqref{secondan} written at the time $s$ by $\dt\sigman(s)$ and integrate over $(0,t)$.
Also in this case, we add the same quantity to both sides.
We thus obtain that
\begin{align}
  & \intQt |\dt\sigman|^2 
  + \frac 12 \, \norma{\sigman(t)}_{B,\tau}^2
  \non
  \\
  & = \frac 12 \, \norma{\sigman(0)}_{B,\tau}^2
  + \intQt \bigl( \sigman - \gammaeps(\phin) \sigman + \kappaeps(\phin)  \bigr) \dt\sigman
  \,\juerg{-}\intQt \Seps \, \phin \, \dt\sigman \,.
  \label{pier2}
\end{align}
Using \eqref{dacauchynbis} for $\sigman$, we can estimate the first term on the \rhs.
Since $\gammaeps$ and $\kappaeps$ are bounded, as well as~$\Seps$, the whole \rhs\ is thus bounded~by
\Beq
  \ceps + \frac 12 \intQt |\dt\sigman|^2
  + \ceps  \intQt |\sigman|^2
  + \ceps \iot \norma{\phin(s)}^2 \, ds \,.
  \non
\Eeq
By accounting for \eqref{primastiman}, and applying the Gronwall lemma, 
we conclude that
\Beq
  \norma\sigman_{\H1H\cap\L\infty{\VB\tau}} \leq \ceps \,.
  \label{secondastiman}
\Eeq

\step
Limit

By \juerg{virtue of} \accorpa{primastiman}{secondastiman} we can find subsequences
(still labeled with the index $n$ for simplicity)
that converge to some limits in the weak or weak star topologies 
\juerg{associated with} the estimates.
Since both $\VA\rho$ and $\VB\tau$ are compactly embedded in~$H$ 
(due to assumption~\eqref{hpAB}), by \pier{recalling, e.g.,
\cite[Sect.~8, Cor.~4]{Simon},}
we conclude that
there is a pair $(\phieps,\sigmaeps)$ satisfying the analogues of the regularity conditions \accorpa{regphi}{regsigma}
such~that\pier{, at least for a subsequence of $n$,}
\Bsist
  & \phin \to \phieps
  & \quad \hbox{weakly star in $\H1H \cap \L\infty{\VA\rho}$}
  \non
  \\
  && \quad \hbox{\pier{and strongly in $\C0H$},} 
  \label{convphin}
  \\
  & \sigman \to \sigmaeps
  & \quad \hbox{weakly star in $\H1H \cap \L\infty{\VB\tau}$}
  \non
  \\
  && \quad \hbox{\pier{and strongly in $\C0H$},} 
  \label{convsigman}
\Esist
We now show that such a pair solves the approximating problem.
Namely, we consider the time integrated version mentioned in Remark~\ref{RemPbleps}.
Take any integer $\bar n\geq1$ and any $v\in\L2{\VAnbar}$.
Then, if the (selected) index~$n$ satisfies $n\geq\bar n$,
then $v(t)\in\VAn$ \aat\ and we can use $v(t)$ as a test function in \eqref{priman}.
After integration over $(0,T)$, we obtain that
\Beq
  \intQ \dt\phin \, v
  + \ioT \bigl( A^\rho\phin(t) , A^\rho v(t) \bigr) \, dt
  + \intQ (\betaeps+\pi)(\phin) \, v
  = \intQ \heps(\phin) \bigl( \meps(\sigman) - \mz \, \ueps \bigr) \, v\,,
  \non
\Eeq
and we can let $n$ tend to infinity.
Since all of the nonlinearities are \Lip\ continuous and $\heps$ and $\meps$,
 as well as~$\ueps$, are even bounded,
we conclude that
\Bsist
  && \intQ \dt\phieps \, v
  + \ioT \bigl( A^\rho\phieps(t) , A^\rho v(t) \bigr) \, dt
  + \intQ (\betaeps+\pi)(\phieps) \, v
  \non
  \\
  && = \intQ \heps(\phieps) \bigl( \meps(\sigmaeps) - \mz \, \ueps \bigr) \, v \,.
  \label{intprimaeps}
\Esist
Since $\bar n$ is arbitrary, this equality holds for every step function $v$
with values in the union $\VAinfty$ of the spaces~$\VAnbar$.
By recalling that $\VAinfty$ is dense in~$\VA\rho$,
an easy density argument shows that \eqref{intprimaeps} holds for every $v\in\L2{\VA\rho}$.

Concerning the equations for $\sigman$ and~$\sigmaeps$,
we observe that $\gammaeps(\phin)$ converges to $\gammaeps(\phi)$ \aeQ, since $\gammaeps$ is continuous.
On the other hand, $\gammaeps$~is bounded.
Hence, we infer that $\gammaeps(\phin)$ also converges to $\gammaeps(\phi)$
in the weak star topology of $\LQ\infty$,
and combining \juerg{this} with the strong convergence given by \eqref{convsigman},
we conclude that $\gammaeps(\phin)\sigman$ converges to $\gammaeps(\phi)\sigmaeps$
weakly in $\L2H$.
Finally, $\kappaeps(\phin)$ converges to \pier{$\kappaeps(\phieps)$} strongly in $\C0H$ by \Lip\ continuity.
By arguing as for the previous equation, and using a similar density property related to the spaces~$\juerg{\VBn}$,
we conclude that $(\phieps,\sigmaeps)$ solves the integrated version of \eqref{secondaeps}
with test functions taken in~$\L2{\VB\tau}$ as well.
Finally, $\phin(0)$ and $\sigman(0)$ converge to $\phieps(0)$ and~$\sigmaeps(0)$ in~$H$ by \accorpa{convphin}{convsigman}.
On the other hand, \eqref{cauchynbis} implies that $\phin(0)$ and $\sigman(0)$ converge to $\phiz$ and~$\sigmaz$,
respectively, in the same topology.
We infer that the initial conditions \eqref{cauchyeps} are satisfied 
and conclude that $(\phieps,\sigmaeps)$ actually solves problem \Pbleps.


\long\def\salta #1\finqui{}

\section{Existence}
\label{EXISTENCE}
\setcounter{equation}{0}

In this section, we prove Theorem~\ref{Existence}.
We start from the solution $(\phieps,\sigmaeps)$ to the approximate problem
and take the limit as $\eps$ tends to zero.
To perform this project, we have to prove some a priori estimates.
In particular, we show that $(\phieps,\sigmaeps)$ satisfies the analogue of \eqref{stimasoluz}
with a constant whose dependence is the same as that of $\overline C_1$ in Theorem~\ref{Existence}
(in~particular, it is independent of~$\eps$),
since the symbol $c$ we use always stands for (possibly different) constants independent of~$\eps$,
according to our general rule.
It follows that this estimate is kept in the limit as $\eps$ goes to zero
and that the last sentence of Theorem~\ref{Existence} is proved as well.
Hence, we do not return to this point in the sequel.

\step 
First a priori estimate

\pier{We claim that from \eqref{pier1}, by a limit procedure as $ n \nearrow \infty$,
it is possible to derive the following inequality}
\Bsist
  && \intQt |\dt\phieps|^2
  + \frac 12 \, \norma{\phieps(t)}_{A,\rho}^2
  + \iO \Betaeps(\phieps(t))
  \non
  \\
  && \pier{\leq{}} \frac 12 \, \norma\phiz_{A,\rho}^2
  + \iO \Betaeps(\phiz)
  + \intQt \bigl(
    \phieps - \pi(\phieps) 
    + \heps(\phieps) \bigl(
      \meps(\sigmaeps) - \mz \, \ueps
    \bigr)
  \bigr) \, \dt\phieps\,.
  \qquad
  \label{testprima}
\Esist
\pier{Indeed, by \eqref{convphin} and the weak lower semicontinuity of norms we infer that
\begin{align*}
&\intQt |\dt\phieps|^2
  + \frac 12 \, \norma{\phieps(t)}_{A,\rho}^2 \leq \liminf_{n \nearrow \infty} 
\intQt |\dt\phin|^2 
  + \liminf_{n \nearrow \infty} \frac 12 \, \norma{\phin(t)}_{A,\rho}^2
\\
&\quad{}
\leq  \liminf_{n \nearrow \infty} \tonde{
\intQt |\dt\phin|^2  +  \frac 12 \, \norma{\phin(t)}_{A,\rho}^2},
\end{align*}
since, $\phin , \phieps$ being weakly continuous from $[0,T]$ to $\VA\rho$, for all $t\in [0,T]$ it occurs that $\phin(t)$ weakly converges to $\phieps(t)$ in  $\VA\rho$. Besides, note that if a sequence $v_n$ converges to $v$ strongly in $H$, then 
$$\iO \Betaeps(v_n ) \to  \iO \Betaeps(v)  \, \hbox{ as } \, n \nearrow \infty.$$
This is due to the mean value theorem in the integral form, which gives 
$$ \iO \Betaeps(v_n ) -   \iO \Betaeps(v)  = \iO \tonde{ \int_0^1 \betaeps (v+ s(v_n - v)) (v_n - v) ds }, $$
and to the \Lip\ continuity of $\betaeps$.  Hence, the terms $\iO \Betaeps(\phin(t))$ 
and  $\iO \Betaeps(\phin(0))$ in \eqref{pier1} converge to the respective ones 
$\iO \Betaeps(\phieps(t))$ 
and  $\iO \Betaeps(\phiz)$ in \eqref{testprima}. Moreover, we point out that
$$ \frac 12 \, \norma{\phin(0)}_{A,\rho}^2 \leq \frac 12 \, \norma{\phiz}_{A,\rho}^2  $$
thanks to \eqref{dacauchynbis}, and that in the last term of \eqref{pier1} we can pass to the limit by strong-weak convergence in $\L2H$ (cf.~\eqref{convphin}--\eqref{convsigman}). Thus, \eqref{testprima} is completely verified.}

All of the terms on the \lhs\ \pier{of \eqref{testprima}} are nonnegative.
As for the ones on the \rhs, we recall \eqref{propBetaeps} for the second one
and use the Schwarz and Young inequality for the volume integral.
Moreover, we notice that 
\eqref{defpsieps} immediately yields that $\psieps$ is uniformly bounded if $\psi$ is bounded,
so that $\meps$ is uniformly bounded since $m$ is so.
Furthermore, by recalling \eqref{defnonlineps}, we account for the inequality \eqref{disugpsitilde} applied to~$\heps$.
Finally, we owe to the inequality $|\ueps|\leq|u|$ \aeQ\ (see \eqref{hpuSeps}).
Then, the \rhs\ of \eqref{testprima} is estimated from above~by
\Bsist
  && c + \frac 12 \intQt |\dt\phieps|^2\, 
	\pier{{}+ \,c\intQt \bigl( 1 + |\phieps|^2 \bigr)
  + c \, \intQt \bigl( 1 + |u |^2 \bigr) \, |\heps(\phieps)|^2}
  \non
  \\
  && \leq \frac 12 \intQt |\dt\phieps|^2
	\juerg{\,+\,c\iot \|\phieps(s)\|^2\,ds}
  + c \, \iot \bigl( 1 + \norma{u(s)}_\infty^2 \bigr) \, \pier{\iO \Betaeps(\phieps(s)) \,ds} + c \,.
  \non
\Esist
Since the function $s\mapsto\norma{u(s)}_\infty^2$ belongs to $L^1(0,T)$ by \eqref{hpuS},
we can apply the Gronwall lemma and conclude that
\Beq
  \norma\phieps_{\H1H\cap\L\infty{\VA\rho}}
  + \norma{\Betaeps(\phieps)}_{\L\infty\Luno}
  \leq c \,.
  \label{primastima}
\Eeq

\step
Second a priori estimate

Similarly, \pier{as for the derivation of \eqref{testprima}, from \eqref{pier2} and \eqref{convphin}--\eqref{convsigman} it is straightforward to deduce the inequality}
\Beq
  \intQt |\dt\sigmaeps|^2
  + \frac 12 \, \norma{\sigmaeps(t)}_{B,\tau}^2
  \pier{{}\leq{} \frac 12 \, \norma\sigmaz_{B,\tau}^2} + \intQt \bigl(
    \sigmaeps
    - \gammaeps(\phieps) \, \sigmaeps
    + \kappaeps(\phieps) 
    - \Seps \phieps
  \bigr) \dt\sigmaeps \,.
  \label{testseconda}
\Eeq
We treat the nontrivial terms on the \rhs, separately.
To do this, we owe to our assumptions \accorpa{hpparam}{hpgrowth}
and to their consequences given by Lemma~\ref{Approx} applied with $\psi=\gamma$ and $\psi=\kappa$.
Moreover, we also account for \eqref{primastima} already established.
We have, by virtue of H\"older's and Young's inequalities, that
\Bsist
  && - \intQt \gammaeps(\phieps) \, \sigmaeps \, \dt\sigmaeps
  \leq \iot \norma{\gammaeps(\phieps(s))}_\pg \, \norma{\sigmaeps(s)}_\qg \, \norma{\dt\sigmaeps(s)}_2 \, ds
  \non
  \\
  && \leq \frac 16 \intQt |\dt\sigmaeps|^2
  + \juerg{c \iot \norma{\gammaeps(\phieps(s))}_\pg^2 \, 
   \norma{\sigmaeps(s)}_\qg^2 \, ds}
  \non
  \\
  && \leq \frac 16 \intQt |\dt\sigmaeps|^2
  + \juerg{\,c\iot \norma{\sigmaeps(s)}_{B,\tau}^2 \, ds} \,,
  \label{Hugo}
\Esist
\juerg{where in the last estimate we have employed 
\eqref{hpembedding} and \eqref{primastima} to
see that}
\Bsist
  &&\juerg{\norma{\gammaeps(\phieps(s))}_\pg^2 
  \,\leq \,c\, \norma{\phieps(s)}_{\ag\pg}^{2\ag} + c
  \,\leq \,c \,\norma{\phieps(s)}_{A,\rho}^{2\ag}+ c
  \,\leq c \quad\mbox{for a.a. $s\in (0,T)$.}} \qquad
  \label{stimagamma}
\Esist
Now, we treat the next term on the right-hand side of \eqref{testseconda}.
By applying the Young inequality, we immediately obtain~that
\Beq
  \intQt \kappaeps(\phieps) \, \dt\sigmaeps
  \leq \frac 16 \intQt |\dt\sigmaeps|^2
  + c \iot \norma{\kappaeps(\phieps(s))}_2^2 \, ds \,.
  \non
\Eeq
On the other hand, we have that
\Bsist
  && \iot \norma{\kappaeps(\phieps(s))}_2^2 \, ds
  \leq c \iot \norma{\phieps(s)}_{2\ak}^{2\ak} \, ds
  + c
  \leq c \iot \norma{\phieps(s)}_{A,\rho}^{2\ak} \, ds
  + c
  \leq c \,. 
  \qquad
  \label{stimakappa}
\Esist
Finally, by using the inequality $|\Seps|\leq|S|$ \aeQ\ (see \eqref{hpuSeps}), we can write
\Beq
  \pier{{}- \intQt \Seps \phieps \, \dt\sigmaeps{}}
  \leq \frac 16 \intQt |\dt\sigmaeps|^2
  + c \intQt |S|^2 \, |\phieps|^2
  \non
\Eeq
and, \pier{due} to assumption \eqref{hpuS} for~$S$ and to \eqref{primastima} once more, 
\Beq
  \intQt |S|^2 \, |\phieps|^2
  \leq \iot \norma{S(s)}_\infty^2 \, \norma{\phieps(s)}^2 \, ds
  \leq c \,.
  \non
\Eeq
By combining all the inequalities just obtained with \eqref{testseconda}
and applying the Gronwall lemma, we conclude~that
\Beq
  \norma\sigmaeps_{\H1H\cap\L\infty{\VB\tau}}
  \leq c \,.
  \label{secondastima}
\Eeq

\step
Estimates of the nonlinear terms

\juerg{By recalling \eqref{stimagamma}, and repeating the arguments that led to 
\eqref{stimakappa} without time integration}, we see on account of \eqref{primastima} that
\Beq
  \norma{\gammaeps(\phieps)}_{\L\infty{\Lx\pg}} 
  \leq c 
  \aand
  \norma{\kappaeps(\phieps)}_{\L\infty H} 
  \leq c \,.
  \label{stimagammakappa}
\Eeq
By combining with \eqref{secondastima} and recalling~\eqref{hpparam},
we deduce that
\Beq
  \norma{\gammaeps(\phieps)\sigmaeps}_{\L\infty H} 
  \leq c \,.
  \label{stimagammasigma}
\Eeq
Since the terms involving~$\heps$ are a little more complicated,
we prepare an auxiliary estimate.
We recall \eqref{hpparam} for the definition of~$\pz$.
Take any $w\in\VA\rho$ and notice that
\Beq
  |\heps(w)|^\ph \leq c \, |w|^{\ah\ph} + c
  \aand
  (|w|\, |\heps(w)|)^\pz
  \leq c \, |w|^{(\ah+1)\pz} + c
  \quad \aeO \,, 
  \non 
\Eeq
thanks to the assumption on $h$ given by~\eqref{hpgrowth}
and Lemma~\ref{Approx} applied with $\psi=h$.
By accounting for \eqref{hpembedding}, we deduce that
\Bsist
  && \norma{\heps(w)}_\ph^\ph
  \leq c \, \norma w_{\ph\ah}^{\ah\ph} + c
  \leq c \, \norma w_{A,\rho}^{\ah\ph} + c
  \aand
  \non
  \\
  && \norma{w\,\heps(w)}_\pz^\pz
  \leq c \, \norma w_{\max\{\ag\pg,2\ak\}}^{\max\{\ag\pg,2\ak\}} + c 
  \leq c \, \norma w_{A,\rho}^{\max\{\ag\pg,2\ak\}} + c \,.
  \non
\Esist
By applying this with $w=\phieps(t)$ \aat\ and recalling \eqref{primastima},
we deduce that
\Beq
  \norma{\heps(\phieps)}_{\L\infty{\Lx\ph}} 
  \leq c 
  \aand
  \norma{\phieps\,\heps(\phieps)}_{\L\infty{\Lx\pz}}
  \leq c \,.
  \non
\Eeq
Since $\meps$ is uniformly bounded (as already observed) 
and $\ueps$ is bounded in $\L2\Linfty$ by \eqref{hpuSeps} and \eqref{hpuS},
we conclude that
\begin{align}
  &\norma{\heps(\phieps)(\meps(\sigmaeps)-\mz\ueps)}_{\L2{\Lx\ph}}
  \non 
  \\
  &{}+ \norma{\heps(\phieps)(\meps(\sigmaeps)-\mz\ueps)\phieps}_{\L2{\Lx\pz}}
  \leq c \,.
  \label{stimah}
\end{align}

\step
Conclusion

We are ready to take the limit of $(\phieps,\sigmaeps)$ as $\eps$ tends to zero.
More precisely, in considering the approximating problem,
we replace \eqref{intprimaeps} (equivalent to~\eqref{primaeps}) and \eqref{secondaeps} 
with the equivalent variational inequality and integrated variational equation, respectively.
We thus start from
\begin{align}
  & \intQ \dt\phieps (\phieps - v)
  + \ioT \bigl( A^\rho\phieps(t) , A^\rho (\phieps(t) - v(t)) \bigr) \, dt
  + \intQ \Betaeps(\phieps)
  + \intQ \pi(\phieps) ( \phieps - v) \
  \non
  \\
  & \leq \intQ \heps(\phieps) \bigl( \meps(\sigmaeps)) - \mz \, \ueps \bigr) (\phieps - v)
  + \intQ \pier{\Betaeps(v)}
  \quad \hbox{for every $v\in\L2{\VA\rho}$}\,,
  \label{disprimaeps}
\end{align}
\pier{and}
\begin{align}  
  & \intQ \dt\sigmaeps \, v
  + \ioT \bigl( B^\tau \sigmaeps(t) , B^\tau v \bigr) \, dt
  + \intQ \gammaeps(\phieps) \, \sigmaeps \, v
  \non
  \\
  & = \intQ \kappaeps(\phieps) \, v
  - \iO \Seps \phieps \, v
  \quad \hbox{for every $v\in\L2{\VB\tau}$}\,,
  \label{intsecondaeps}
\end{align}
\pier{as well as} the initial conditions~\eqref{cauchyeps}.
By \eqref{primastima} and \eqref{secondastima}, 
and accounting for standard weak, weak star and strong compactness results
(see, e.g., \cite[Sect.~8, Cor.~4]{Simon} for the latter),
we see that there exists a strictly decreasing subsequence~$\eps_n\searrow0$ such that the corresponding convergence holds true.
However, we still write~$\eps$, at least for a while, to simplify the notation.
Namely, we have that
\Bsist
  & \phieps \to \phi
  & \quad \hbox{weakly star in $\H1H\cap\L\infty{\VA\rho}$},
  \non
  \\
  && \qquad \hbox{strongly in $\C0H$ and \aeQ},
  \label{convphi}
  \\
  & \sigmaeps \to \sigma
  & \quad \hbox{weakly star in $\H1H\cap\L\infty{\VB\tau}$},
  \non
  \\
  && \qquad \hbox{strongly in $\C0H$ and \aeQ}, 
  \label{convsigma}
\Esist
for some pair $(\phi,\sigma)$ satisfying \accorpa{regphi}{regsigma}.
It follows that such a pair fulfills the initial conditions \eqref{cauchy}
and that
\Beq
  \pi(\phieps) \to \pi(\phi)
  \quad \hbox{strongly in $\C0H$},
  \label{convpi}
\Eeq
because $\pi$ is \Lip\ continuous.
Moreover, since $\pg$, $\ph$ and $\pz$ are larger than~$1$, $\meps$ is uniformly bounded, 
and the estimates \accorpa{stimagammakappa}{stimah} hold true,
we also have that
\Bsist
  & \gammaeps(\phieps)\sigmaeps \to \zeta_1
  & \quad \hbox{weakly star in \juerg{$L^\infty(0,T;H)$}},
  \label{convg}
  \\
  & \kappaeps(\phieps) \to \zeta_2
  & \quad \hbox{weakly star in $\L\infty H$},
  \label{convk}
  \\
  & \heps(\phieps)(\meps(\sigmaeps)-\mz\ueps) \to \zeta_3
  & \quad \hbox{\pier{weakly in $\L2{\Lx\ph}$}},
  \label{convhphi}
  \\
  & \heps(\phieps)(\meps(\sigmaeps)-\mz\ueps)\phieps \to \zeta_4
  & \quad \hbox{\pier{weakly in $\L2{\Lx\pz}$}} ,
  \label{convhphibis}
  \\
  & \meps(\sigmaeps) \to \zeta_5
  & \quad \hbox{weakly star in $\LQ\infty$},
  \label{convm}
\Esist
for a suitable subsequence and some limiting functions~$\zeta_i$.
On the other hand, $\gammaeps$ converges to $\gamma$ uniformly on every compact interval of~$\erre$
thanks to Lemma~\ref{Approx} applied with $\psi=\gamma$.
Combining this with the almost everywhere convergence given by \accorpa{convphi}{convsigma},
we deduce that
\Beq
  \gammaeps(\phieps)\sigmaeps \to \gamma(\phi) \sigma
  \quad \aeQ \,.
  \non
\Eeq
We infer that \pier{(see, e.g., \cite[Lemme~1.3, p.~12]{Lions})} $\, \zeta_1=\gamma(\phi)\sigma \, $ \aeQ.
By analogously arguing for the other nonlinear terms appearing in \accorpa{convk}{convm}
(and recalling \eqref{hpuSeps} for~$\ueps$),
we also conclude that
\Beq
  \zeta_2 = \kappa(\phi) \,, \quad
  \zeta_3 = h(\phi) \bigl( m(\sigma) - \mz u \bigr) \,, \quad
  \zeta_4 = \zeta_3 \, \phi
  \aand
  \zeta_5 = m(\phi) \,,
  \non
\Eeq
and that \juerg{these} limits can also be understood in the sense of convergence \aeQ.
Furthermore, as $\pz>1$, \juerg{the convergence in} \eqref{convhphibis} 
also holds in the strong topology of~$\LQ1$.
It is also clear that \eqref{hpuSeps}, \eqref{hpuS} and \eqref{convphi} 
imply that $\Seps\phieps$ converges to~$S\phi$, e.g., weakly in $\L2H$ 
(by~boundedness in this space and convergence \aeQ).
\pier{For the last term of \eqref{disprimaeps} we note that 
$$
\intQ \pier{\Betaeps(v)} \leq \intQ \pier{\Beta(v)}
\quad \hbox{ for every $v\in\L2{\VA\rho}$}
$$
by \eqref{propBetaeps}. Then, all of the above} ensures that we can take the limit 
in all of the terms of \accorpa{disprimaeps}{intsecondaeps}
but the second and third ones on the \lhs\ of the variational inequality.
In particular, we obtain the time-integrated version of~\eqref{seconda}
with test functions $v$ in $\L2{\VB\tau}$.
For the first of the terms of \eqref{disprimaeps} we still have to consider, we have that
\Beq
  \ioT \bigl( A^\rho\phi(t) , A^\rho (\phi(t) - v(t)) \bigr) \, dt
  \leq \liminf_{\eps\searrow0} \ioT \bigl( A^\rho\phieps(t) , A^\rho (\phieps(t) - v(t)) \bigr) \, dt\,,
  \label{semicont}
\Eeq
by the semicontinuity of the norm and the weak convergence of~$A^\rho\phieps$ to $A^\rho\phi$
in $\L2H$ ensured by \eqref{convphi}.
For the other one, we are going to derive the inequality
\Beq
  \iO \Beta(\phi(t)) 
  \leq \liminf_{\eps\searrow0} \iO \Betaeps(\phieps(t)) 
  \quad \aat \,.
  \label{liminfBeta}
\Eeq
Notice that its \rhs\ (as a function of~$t$) is bounded by~\eqref{primastima}.
In particular, the requirement $\Beta(\phi)\in\L\infty\Luno$ (see~\eqref{regBetaphi}) is fulfilled
once \eqref{liminfBeta} is established,
and we also have that
\Beq
  \intQ \Beta(\phi) 
  \leq \liminf_{\eps\searrow0} \intQ \Betaeps(\phieps)\,, 
  \label{liminfBetaQ}
\Eeq
by Fatou's lemma applied to the functions $t\mapsto\iO\Betaeps(\phieps(t))$.
It is clear that \accorpa{semicont}{liminfBetaQ}
are understood for the subsequence $\{\eps_n\}$ selected in relation to \accorpa{convphi}{convm}.
In proving \eqref{liminfBeta}, we use this subsequence in the notation, for clarity.
Moreover, we account for the properties \eqref{propBetaeps} of~$\Betaeps$.
Let us start.
Since $\epsn<\epsm$ whenever $n>m$, we have~that
\Beq
  \Betaepsm(\phiepsn(t))
  \leq \Betaepsn(\phiepsn(t))
  \quad \hbox{\aeQ, for every $n>m$},
  \non
\Eeq
whence also 
\Beq
  \Betaepsm(\phi)
  = \lim_{n\to\infty} \Betaepsm(\phiepsn)
  = \liminf_{n\to\infty} \Betaepsm(\phiepsn)
  \leq \liminf_{n\to\infty} \Betaepsn(\phiepsn)
  \quad \aeQ ,
  \non
\Eeq
since $\Betaepsm$ is continuous.
On the other hand, we have that
\Beq
  \Beta(\phi) = \lim_{m\to\infty} \Betaepsm(\phi)
  \quad \aeQ.
\Eeq
We infer that
\Beq
  \Beta(\phi(t))
  \leq \liminf_{n\to\infty} \Betaepsn(\phiepsn(t))
  \quad \aeO ,\quad \aat ,
\Eeq
and \eqref{liminfBeta} follows from Fatou's lemma.
At this point, we can owe to \eqref{liminfBetaQ} and let $\eps$ tend to zero in \eqref{disprimaeps} as well
to obtain \eqref{intprima} with test functions $v$ in $\L2{\VA\rho}$.
We conclude that $(\phi,\sigma)$ actually solves problem \Pbl ,
and the proof of Theorem~\ref{Existence} is complete.


\section*{Acknowledgments}
\pier{This research was supported by the Italian Ministry of Education, 
University and Research~(MIUR): Dipartimenti di Eccellenza Program (2018--2022) 
-- Dept.~of Mathematics ``F.~Casorati'', University of Pavia. 
In addition, {PC gratefully mentions} some other support 
from the GNAMPA (Gruppo Nazionale per l'Analisi Matematica, 
la Probabilit\`a e le loro Applicazioni) of INdAM (Isti\-tuto 
Nazionale di Alta Matematica).}


\vspace{3truemm}

\Begin{thebibliography}{10}

{\small%
\bibitem{Baio1}		
C. Baiocchi, Sulle equazioni differenziali astratte lineari del primo e del secondo ordine negli spazi di Hilbert, Ann. Mat. Pura Appl. (4)  {\bf 76}  (1967), 233--304.

\bibitem{Baio2}
C. Baiocchi, Soluzioni ordinarie e generalizzate del problema di Cauchy per equazioni differenziali astratte lineari del secondo ordine in spazi di Hilbert,
Ricerche Mat. {\bf 16}  (1967), 27--95.

{\bibitem{BKM}
B. Baeumer, M. Kov\'acs, M. M. Meerschaert, 
Numerical solutions for fractional reaction-diffusion equations,
Comput. Math. Appl.  {\bf 55}  (2008), 2212--2226.}

\bibitem{Barbu}
V. Barbu,
``Nonlinear Differential Equations of Monotone Type in Banach Spaces'',
Springer, London, New York, 2010.

\bibitem{BLM}
N. Bellomo, N. K. Li, P. K. Maini,
On the foundations of cancer modelling: Selected topics, speculations, and perspectives,
{Math. Models Methods Appl. Sci.} \textbf{18} (2008), 593--646.

\bibitem{Brezis}
H. Brezis,
``Op\'erateurs maximaux monotones et semi-groupes de contractions
dans les espaces de Hilbert'',
North-Holland Math. Stud. {\bf 5},
North-Holland,
Amsterdam,
1973.

\bibitem{CRW}
C.~Cavaterra, E.~Rocca, H.~Wu, {Long-time dynamics and optimal
  control of a diffuse interface model for tumor growth}, {Appl. Math.
  Optim.}, DOI: 10.1007/s00245-019-09562-5.

{\bibitem{CB}
S. K. Chandra, M. K. Bajpai, 
Mesh free alternate directional implicit method based three dimensional 
super-diffusive model for benign brain tumor segmentation,
{Comput. Math. Appl.}  {\bf 77}  (2019), 3212--3223.}

\bibitem{CWSL}
Y. Chen, S. M. Wise, V. B. Shenoy, J. S. Lowengrub,
A stable scheme for a nonlinear multiphase tumor growth model with an elastic membrane,
{Int. J. Numer. Methods Biomed. Eng.} \textbf{30} (2014), 726--754.

\bibitem{CG1} 
P. Colli, G. Gilardi, 
Well-posedness, regularity and asymptotic analyses 
for a fractional phase field system,
{Asymptot. Anal.} \textbf{114} (2019), 93--128.

\bibitem{CGH}
P.~Colli, G.~Gilardi, D.~Hilhorst, {On a {Cahn-Hilliard} type phase
  field system related to tumor growth}, {Discrete Contin. Dyn. Syst. Ser.
  A} \textbf{35} (2015), 2423--2442.

\bibitem{CGMR3}
P.~Colli, G.~Gilardi, G.~Marinoschi, E.~Rocca, {Sliding mode control
  for a phase field system related to tumor growth}, {Appl. Math. Optim.}
  \textbf{79} (2019), 647--670.

\bibitem{CGRS3}
P. Colli, G. Gilardi, E. Rocca, J. Sprekels,
Optimal distributed control of a diffuse interface model of tumor growth,
{Nonlinearity} \textbf{30} (2017), 2518--2546.

\bibitem{CGS18}
P. Colli, G. Gilardi, J. Sprekels,
Well-posedness and regularity for a generalized fractional Cahn--Hilliard system,
{Atti Accad. Naz. Lincei Rend. Lincei Mat. Appl.} {\bf 30} (2019), 437--478.

\bibitem{CGS23}
P. Colli, G. Gilardi, J. Sprekels,
Well-posedness and regularity for a fractional tumor growth model, 
{Adv. Math. Sci. Appl.} {\bf 28} (2019), 343--375.

\bibitem{CGS25}
P. Colli, G. Gilardi, J. Sprekels,
A distributed control problem for a fractional tumor growth model,
Mathematics {\bf 7} (2019), 792.

\bibitem{CGS21bis}
P. Colli, G. Gilardi, J. Sprekels,
Recent results on well-posedness and optimal control 
for a class of generalized fractional Cahn--Hilliard
systems, Control Cybernet. {\bf 48} (2019), 153--197. 

\bibitem{CGS19}
P. Colli, G. Gilardi, J. Sprekels,
Optimal distributed control of a generalized 
fractional Cahn--Hilliard system,
Appl. Math. Optim. {\bf 82} (2020), 551--589.

\bibitem{CGS22}
P. Colli, G. Gilardi, J. Sprekels,
Longtime behavior for a generalized Cahn--Hilliard 
system with fractional operators,
Atti Accad. Peloritana Pericolanti Cl. Sci. Fis. Mat. 
Natur. {\bf 98} (2020), suppl.~2, A4, 18 pp.

\bibitem{CGS24}
P. Colli, G. Gilardi, J. Sprekels, 
{Asymptotic analysis of a tumor growth model 
with fractional operators}, 
Asymptot. Anal. {\bf 120} (2020), 41--72.

\bibitem{CGS21}
P. Colli, G. Gilardi, J. Sprekels,
Deep quench approximation and optimal control of 
general Cahn--Hilliard systems with fractional
operators and double-obstacle potentials, 
Discrete Contin. Dyn. Syst. Ser. S {\bf 14} (2021), 243--271.

\bibitem{CGS27}
P. Colli, G. Gilardi, J. Sprekels, 
{Optimal control of a phase field system of Caginalp type with 
fractional operators}, Pure Appl. Funct. Anal., to appear
(see also preprint  arXiv:2005.11948 [math.AP] (2020), pp.~1--38).

\bibitem{CGLMRR1}
P. Colli, H. Gomez, G. Lorenzo, G. Marinoschi, A. Reali, E. Rocca,
Mathematical analysis and simulation study of a phase-field model 
of prostate cancer growth with chemotherapy and antiangiogenic therapy effects,
Math. Models Methods Appl. Sci. {\bf 30} (2020), 1253--1295.

\bibitem{CGLMRR2}
P. Colli, H. Gomez, G. Lorenzo, G. Marinoschi, A. Reali, E. Rocca, 
{Optimal control of cytotoxic and antiangiogenic 
therapies on prostate cancer growth},  Math. Models Methods Appl. 
Sci., to appear (see also preprint arXiv:2007.05098 [math.OC] (2020), pp.~1--41).

\bibitem{CL2010}
V. Cristini, J. S. Lowengrub,
``Multiscale Modeling of Cancer: An Integrated Experimental 
and Mathematical Modeling Approach'',
Cambridge Univ. Press, Cambridge, 2010.

\bibitem{CLLW}
V. Cristini, X. Li, J. S. Lowengrub, S. M. Wise, 
Nonlinear simulations of solid tumor growth using a 
mixture model: invasion and branching, 
{J. Math. Biol.} {\bf 58} (2009), 723--763.  

\bibitem{DFRSS}
M. Dai, E. Feireisl, E. Rocca, G. Schimperna, M. Schonbek,
Analysis of a diffuse interface model for multi-species tumor growth,
{Nonlinearity} \textbf{30} (2017), 1639--1658.

\bibitem{EGAR}
M. Ebenbeck, H. Garcke,
Analysis of a Cahn--Hilliard--Brinkman model for tumour growth with chemotaxis,
{J. Differential Equations} \textbf{266} (2019), 5998--6036.

\bibitem{EGPS}
G. Estrada-Rodriguez, H. Gimperlein, K. J. Painter, J. Stocek,  
Space-time fractional diffusion in cell movement models with delay.
{Math. Models Methods Appl. Sci.}  {\bf 29} (2019), 65--88.

\bibitem{EL}
L. R. Evangelista, E. K. Lenzi, 
``Fractional diffusion equations and anomalous diffusion'',
Cambridge University Press, Cambridge, 2018.

\bibitem{FBG2006}
A. Fasano, A. Bertuzzi, A. Gandolfi,
Mathematical modeling of tumour growth and treatment,
Complex Systems in Biomedicine, Springer, Milan, 2006, {pp.~71--108.}

\bibitem{Frieboes2010}
H.~B. Frieboes, F.~Jin, Y.-L. Chuang, S.~M. Wise, J.~S. Lowengrub and
  V.~Cristini, {Three-dimensional multispecies nonlinear tumor
  growth--{II}: Tumor invasion and angiogenesis}, {J. Theor. Biol.}
  \textbf{264} (2010), 1254--1278.

\bibitem{Fri2007}
A. Friedman,
Mathematical analysis and challenges arising from models of tumor growth,
{Math. Models Methods Appl. Sci.} \textbf{17} (2007), 1751--1772.

\bibitem{FGR}
S.~Frigeri, M.~Grasselli and E.~Rocca, {On a diffuse interface model of
  tumour growth}, {European J. Appl. Math.} \textbf{26} (2015), 215--243.

\bibitem{FLR}
S. Frigeri, K.~F. Lam, E. Rocca, 
On a diffuse interface model for tumour growth with non-local 
interactions and degenerate mobilities, 
in ``Solvability, regularity, and optimal control of boundary 
value problems for PDEs'', 
P.~Colli, A.~Favini, E.~Rocca, G.~Schimperna, J.~Sprekels~(ed.), 
Springer INdAM Series~{\bf 22}, Springer, Cham, 2017, pp.~217--254.

\bibitem{GLR}
H. Garcke, K.~F. Lam, E. Rocca,
Optimal control of treatment time in a diffuse interface model for tumour growth,
{Appl. Math. Optim.} {\bf 78} (2018), 495--544.

\bibitem{GLS}
H. Garcke, K. F. Lam, A. Signori,
On a phase field model of Cahn--Hilliard type for tumour growth with mechanical effects,
Nonlinear Anal. Real World Appl. {\bf 57} (2021), 103192, 28 pp.

\bibitem{G-B}
R. Granero-Belinch\'on, 
Global solutions for a hyperbolic-parabolic system of chemotaxis,
{J. Math. Anal. Appl.}  {\bf 449} (2017), 872--883.

\bibitem{HDPZO}
A. Hawkins-Daarud, S. Prudhomme, K.~G. van der Zee, J.~T. Oden,
Bayesian calibration, validation, and uncertainty quantification of diffuse 
interface models of tumor growth,
{J. Math. Biol.} {\bf 67} (2013), 1457--1485.

\bibitem{HZO12}
A. Hawkins-Daarud, K. G. van der Zee, J. T. Oden,
Numerical simulation of a thermodynamically consistent four-species tumor growth model,
{Int. J. Numer. Meth. Biomed. Engrg.} \textbf{28} (2012), 3--24.

\bibitem{INK}
R. W. Ibrahim, H. K. Nashine, N. Kamaruddin, 
Hybrid time-space dynamical systems of growth bacteria with applications in segmentation,
{Math. Biosci.}  {\bf 292}  (2017), 10--17.

\bibitem{JWZ}
J. Jiang, H. Wu, S. Zheng,
Well-posedness and long-time behavior of a non-autonomous
Cahn--Hilliard--Darcy system with mass source modeling tumor growth,
J. Differential Equations {\bf 259} (2015), 3032--3077.

\bibitem{KU}		
K. H. Karlsen, S. Ulusoy, 
On a hyperbolic Keller--Segel system with degenerate nonlinear fractional diffusion,
{Netw. Heterog. Media}  {\bf 11} (2016), 181--201.

\bibitem{Lima2014}
E.~A. B.~F. Lima, J.~T. Oden and R.~C. Almeida, {A hybrid ten-species
  phase-field model of tumor growth}, {Math. Models Methods Appl. Sci.}
  \textbf{24} (2014) 2569--2599.

\bibitem{Lions}
J.-L. Lions, ``Quelques M\'ethodes de R\'esolution 
des Probl\`emes aux Limites non Lin\'eaires'', 
Dunod Gauthier-Villars, Paris, 1969. 

\bibitem{lorenzo2016tissue}
G.~Lorenzo, M.~A. Scott, K.~Tew, T.~J.~R. Hughes, Y.~J. Zhang, L.~Liu,
  G.~Vilanova and H.~Gomez, {Tissue-scale, personalized modeling and
  simulation of prostate cancer growth}, {Proc. Natl. Acad. Sci. U.S.A.}
  \textbf{113} (2016) E7663--E7671.

\bibitem{LTZ}
J. S. Lowengrub, E. S. Titi, K. Zhao,
Analysis of a mixture model of tumor growth,
{European J. Appl. Math.} \textbf{24} (2013), 691--734.

\bibitem{MRS}
A.~Miranville, E.~Rocca and G.~Schimperna, {On the long time behavior
  of a tumor growth model}, {J. Differential Equations} \textbf{267} (2019),
  2616--2642.

\bibitem{OHP}
J.~T. Oden, A. Hawkins, S. Prudhomme,
General diffuse-interface theories and an approach to predictive tumor growth modeling,
{Math. Models Methods Appl. Sci.} {\bf 20} (2010), 477--517. 

\bibitem{Sig}
A. Signori,
Optimal distributed control of an extended model of tumor 
growth with logarithmic potential,
Appl. Math. Optim.  {\bf 82} (2020), 517--549.

\bibitem{Simon}
J. Simon,
Compact sets in the space $L^p(0,T; B)$,
Ann. Mat. Pura Appl.~(4)
{\bf 146} (1987), 65--96.

\bibitem{SAJM}
A. Sohail, S. Arshad, S. Javed, K. Maqbool, 
Numerical analysis of fractional-order tumor model,
{Int. J. Biomath.}  {\bf 8} (2015), 1550069, 13 pp.

\bibitem{SA}
N. H. Sweilam, S. M. Al-Mekhlafi,   
Optimal control for a nonlinear mathematical model of tumor under 
immune suppression: a numerical approach,
{Optimal Control Appl. Methods} {\bf 39} (2018), 1581--1596.

\bibitem{Wise2011}
S. M. Wise, J. S. Lowengrub, V. Cristini,
An adaptive multigrid algorithm for simulating solid tumor growth using mixture models,
{Math. Comput. Modelling} \textbf{53} (2011), 1--20.

\bibitem{Wise2008}
S. M. Wise, J. S. Lowengrub, H. B. Frieboes, V. Cristini,
Three-dimensional multispecies nonlinear tumor growth - I: model and numerical method,
{J. Theoret. Biol.} {\bf 253} (2008), 524--543.

\bibitem{Xu2016}
J.~Xu, G.~Vilanova and H.~Gomez, {A mathematical model coupling tumor
  growth and angiogenesis}, {PLoS ONE} \textbf{11} (2016) e0149422.

{\bibitem{ZSMD}
Y. Zhou, L. Shangerganesh, J. Manimaran, A. Debbouche, 
A class of time-fractional reaction-diffusion equation with nonlocal boundary condition,
{Math. Methods Appl. Sci.}  {\bf 41} (2018), 2987--2999.}

}

\End{thebibliography}

\End{document}
